\documentclass[12pt,twoside,reqno]{amsart}
\usepackage{xcolor}
\usepackage{graphicx}
\usepackage{mathtools}
\usepackage{mathrsfs}
\usepackage{cite}
\usepackage{cleveref}
\usepackage{geometry}
\usepackage{marginnote}
\usepackage{todonotes}
\newcommand{\blm}[1]{{#1}}
\newtheorem{theorem}{Theorem}[section]
\newtheorem{corollary}[theorem]{Corollary}
\newtheorem{lemma}[theorem]{Lemma}
\newtheorem{proposition}[theorem]{Proposition}

\allowdisplaybreaks
\numberwithin{equation}{section}
\begin{document}
\title{Classification of radially symmetric solutions to\\ a Hardy-H\'enon equation with sublinear absorption\\ and applications}
\thanks{}
\author{Razvan Gabriel Iagar}
\address{Departamento de Matem\'{a}tica Aplicada, Ciencia e Ingenieria de los Materiales y Tecnologia Electr\'onica, Universidad Rey Juan Carlos, M\'{o}stoles, 28933, Madrid, Spain}
\email{razvan.iagar@urjc.es}
\author{Philippe Lauren\c{c}ot}
\address{Laboratoire de Math\'ematiques (LAMA) UMR~5127, Universit\'e Savoie Mont Blanc, CNRS, F--73000 Chamb\'ery, France}
\email{philippe.laurencot@univ-smb.fr}


\date{\today}

\begin{abstract}
Uniqueness of a non-negative radially symmetric solution with non-increasing profile to the Hardy-H\'enon elliptic equation with sublinear absorption
\begin{equation*}
	\Delta v(x)+|x|^{\sigma}v(x)-v^{1/m}(x)=0, \quad x\in\mathbb{R}^N,
\end{equation*}
is established in the range $N\geq 1$, $m>1$ and $\sigma\in \big(\max\{-N,-2\},0\big)$. Consequences of this result are the uniqueness (up to scaling) of the non-negative extremal function to a Caffarelli-Kohn-Nirenberg inequality and the uniqueness of a non-negative radially symmetric self-similar solution with non-increasing profile to a parabolic Hardy-H\'enon equation with quasilinear diffusion. The proof relies on ordinary differential equations techniques, which can be used thanks to the radially symmetric setting, and requires in particular a delicate analysis of the well-posedness and continuous dependence in view of the non-Lipschitz nonlinearity and the singularity of $|x|^\sigma$ at the origin.
\end{abstract}

\maketitle

%
%
\pagestyle{myheadings}
\markboth{\sc{R.G. Iagar \& Ph. Lauren\c cot}}{\sc{Uniqueness of solutions to a Hardy-H\'enon equation}}

\section{Introduction}

The purpose of this work is to show the uniqueness of a non-negative radially symmetric solution $v\in \mathcal{X} := H^1(\mathbb{R}^N)\cap L^{(m+1)/m}(\mathbb{R}^N)$ with non-increasing profile to the Hardy-H\'enon elliptic equation with sublinear absorption 
\begin{equation}\label{eq1}
	\Delta v(x) + |x|^{\sigma} v(x) - v^{1/m}(x) = 0, \quad x\in\mathbb{R}^N,
\end{equation}
in the range of dimension and exponents
\begin{equation}\label{range.exp}
	N\ge 1, \quad m>1, \quad \sigma\in \big(\max\{-N,-2\},0\big).
\end{equation}
Existence of such solutions is established in \cite{IL25}, where it is additionally shown that they belong to $L^\infty(\mathbb{R}^N)$ and are compactly supported. The motivation to study non-negative radially symmetric solutions to~\eqref{eq1} in $\mathcal{X}$ with non-increasing profile is twofold, as described below.
	
\medskip
	
\paragraph{\textbf{Caffarelli-Kohn-Nirenberg (CKN) inequalities}} Setting $p:=(m+1)/m\in (1,2)$ and introducing the functional
\begin{equation}
	\mathcal{S}(w) := \frac{\|\nabla w\|_2^{2a} \|w\|_{p}^{2(1-a)}}{N_\sigma^2(w)}, \quad w\in\mathcal{X}, \ w\ne 0, \label{ckn01}
\end{equation}
where
\begin{equation}
	a := \frac{N(m-1)-\sigma(m+1)}{N(m-1)+2(m+1)} = \frac{N(2-p)-\sigma p}{N(2-p)+2p}\in (0,1), \label{ckne}
\end{equation}
and
\begin{equation*}	
	N_\sigma^2(w) := \int_{\mathbb{R}^N} |x|^\sigma w^2(x)\,dx,
\end{equation*}
it follows from \cite[Theorem]{CKN84} that $\mathcal{X}$ embeds continuously in $L^2(\mathbb{R}^N,|x|^\sigma\,dx)$, thereby guaranteeing the finiteness of $N_\sigma(w)$ for all $w\in\mathcal{X}$, and the positivity and finiteness of
\begin{equation}
	\mathcal{S}_* := \inf_{w\in\mathcal{X}, w\ne 0} \mathcal{S}(w)\in (0,\infty). \label{ckni}
\end{equation}
According to \cite{IL25,ZZ15}, there are minimizers $v\in \mathcal{X}$ to~\eqref{ckni} satisfying $\mathcal{S}(v) = \mathcal{S}_*$, also referred to as extremal functions, and all non-negative minimizers of~\eqref{ckni} are non-negative radially symmetric functions in $\mathcal{X}$ with non-increasing profile, as it follows from Proposition~\ref{prop.ckn01}. Moreover, given such a non-negative minimizer $v$, we show in Proposition~\ref{prop.ckn01} that there are $(\zeta,\eta)\in (0,\infty)^2$ such that $x\mapsto \zeta v(\eta x)$ is a solution to~\eqref{eq1} in $\mathcal{X}$. Denoting the unique non-negative radially symmetric solution to~\eqref{eq1} in $\mathcal{X}$ with non-increasing profile by $V_*$, a consequence of the uniqueness result established in Theorem~\ref{th.uniq} is that the set of non-negative minimizers of~\eqref{ckni} is 
\begin{equation*}
	\mathcal{M} := \left\{ x\mapsto \zeta V_*(\eta x), \ (\zeta,\eta)\in (0,\infty)^2 \right\}.
\end{equation*}	
We deduce that, up to scaling invariance, there is a unique extremal function for the CKN inequality~\eqref{ckni}. Additional properties derived in \cite{IL25} and Theorem~\ref{th.uniq} are the boundedness of $V_*$ and the compactness of its support, along with its behavior at the origin and at the edge of its support. It is worth mentioning at this point that, to the best of our knowledge, uniqueness (up to scaling invariance) of extremal functions for CKN inequalities seems to be only available in a few cases, where the minimizers are actually explicit, see \cite{CW01,CC09,CC22,DEL16,DELM17,DFLLxx} and the references therein.

\medskip
	
\paragraph{\textbf{A quasilinear Hardy-H\'enon parabolic equation with source}} The study of~\eqref{eq1} is also motivated by the quest for understanding the dynamics of non-negative solutions to the parabolic Hardy-H\'enon equation with quasilinear diffusion
\begin{equation}\label{eq.parab}
	\partial_t u = \Delta u^m + |x|^{\sigma} u^m, \quad t>0,\ x\in \mathbb{R}^N,
\end{equation}
with $(m,\sigma)$ satisfying~\eqref{range.exp}. We proved in \cite[Proposition~4.3]{IL26} that non-trivial, non-negative solutions to~\eqref{eq.parab} blow up in a finite time. In view of the already known results dealing with the spatially homogeneous case $\sigma=0$ \cite{S4}, the blow-up is expected to be self-similar of type~I and all (non-negative) self-similar solutions to~\eqref{eq.parab} (if any) are in separate variables backward form
\begin{equation}\label{SSS}
	u(t,x) = (T-(m-1)t)^{-1/(m-1)} U(|x|), \quad (t,x)\in \left[0,\frac{T}{m-1} \right) \times\mathbb{R}^N,
\end{equation}
for some arbitrary $T\in(0,\infty)$, as indicated in \cite{IL25}. Substituting the ansatz~\eqref{SSS} into~\eqref{eq.parab}, easy algebraic manipulations show that the profile $U^m$ is a non-negative radially symmetric solution to~\eqref{eq1}. Therefore, under the additional requirements that $x\mapsto U^m(|x|)$ belongs to $\mathcal{X}$ and $U^m$ is non-increasing, the uniqueness result established in Theorem~\ref{th.uniq} implies that $U^m$ coincides with the unique non-negative radially symmetric solution to~\eqref{eq1} in $\mathcal{X}$ with non-increasing profile $V_*$. Such a result paves the way towards a complete description of the behavior of non-negative solutions to~\eqref{eq.parab} at the blow-up time, at least for a suitably chosen class of initial conditions.

\medskip

We thus study non-negative radially symmetric solutions $v$ to~\eqref{eq1}. Letting $r=|x|$ and $V(r)=v(|x|)$, $x\in\mathbb{R}^N$,  it follows from~\eqref{eq1} that the profile $V$ solves
\begin{equation}\label{eq0.rad}
	V''(r)+\frac{N-1}{r}V'(r)+r^{\sigma}V(r)-V^{1/m}(r)=0, \quad r>0,
\end{equation}
and our aim is to prove that~\eqref{eq0.rad} has a unique non-negative and non-increasing solution on $(0,\infty)$. To this end, we employ a shooting technique and classify the solutions to~\eqref{eq0.rad} with initial conditions
\begin{equation}
	V(0)=a>0, \ \lim_{r\to 0} r^{N-1} V'(r) = 0, \label{ica0}
\end{equation}
according to the value of $a>0$. We actually face a first difficulty in the analysis of the well-posedness of~\eqref{eq0.rad}, as the singularity of $r\mapsto r^\sigma$ as $r\to 0$ prevents us from applying the classical Cauchy-Lipschitz theorem. Establishing the well-posedness of~\eqref{eq0.rad}, supplemented with the initial condition~\eqref{ica0}, is one step of the proof and requires the unusual condition on the behavior of $V'(r)$ as $r\to 0$. In addition, Eq.~\eqref{eq0.rad} features a non-Lipschitz nonlinearity due to $m>1$, which also prevents us from using the classical Cauchy-Lipschitz theorem where $V$ vanishes, a phenomenon which occurs for a whole interval of values of $a$. In particular, continuous dependence with respect to $a$ is not granted by classical results and \blm{needs} to be proved prior to its use. Having dealt with these issues, we next prove that, for all $a>0$, there is $R(a)\in (0,\infty)$ such that
\begin{equation}
	V(r)>0 \;\text{ and }\; V'(r)<0, \quad r\in (0,R(a)), \quad\qquad V(R(a))V'(R(a))=0. \label{pim}
\end{equation}
At this point, the behavior of $V(r)$ as $r\to R(a)$ depends on $a$ and three different cases are to be \blm{handled}:
\begin{itemize}
	\item [(a)] $V(R(a))=0$ and $V'(R(a))<0$. In that case, $x\mapsto V(|x|)$ cannot be extended to a non-negative solution to~\eqref{eq1} in $\mathbb{R}^N$, though it might be extended to a sign-changing solution to the elliptic equation $\Delta v(x) + |x|^{\sigma} v(x) - \big(|v|^{(1-m)/m}v\big)(x) = 0$, $x\in\mathbb{R}^N$;
	\item [(b)] $V(R(a))=V'(R(a))=0$. The existence of such a solution is established in \cite{IL25} and, once extended by zero to $(R(a),\infty)$, it is the profile of the sought-for non-negative radially symmetric solution to~\eqref{eq1} with non-increasing profile. The main contribution of this paper is that the property $V(R(a))=V'(R(a))=0$ holds true for a single value of the shooting parameter $a$;
	\item [(c)] $V(R(a))>0$ and $V'(R(a))=0$. In that case, we prove that $V'$ is positive in a right neighborhood \blm{of $R(a)$}, so that $x\mapsto V(|x|)$ cannot be a non-negative radially symmetric solution to~\eqref{eq1} with non-increasing profile.
\end{itemize}

We are now in a position to make precise our main results.

\section{Main results}

As previously mentioned, our goal is to classify non-negative solutions to
\begin{subequations}\label{CPV} 
\begin{equation}\label{eq1.rad}
	V''(r)+\frac{N-1}{r}V'(r)+r^{\sigma}V(r)-V^{1/m}(r)=0, \quad r>0,
\end{equation}
supplemented with the following initial conditions
\begin{equation}\label{ica}
	V(0)=a>0, \quad \lim\limits_{r\to 0}r^{N-1}V'(r)=0\blm{,}
\end{equation}
\end{subequations}
and establish the uniqueness of a non-negative and non-increasing solution to~\eqref{CPV} presenting an interface.  We observe that the second condition in~\eqref{ica} is not a standard one. Indeed, despite the radial symmetry, one cannot assume $V'(0)=0$ as such a solution might not exist, in particular for $\sigma\in(-2,-1)$, as established in \cite{IL25}. This possibly singular behavior as $r\to 0$ of solutions to the Cauchy problem~\eqref{CPV} is an effect of the negativity of $\sigma$ (and thus of the singular character of the coefficient $r^{\sigma}$ as $r\to 0$). However, the condition \eqref{ica} ensures the well-posedness of the problem. We are now ready to state the main results of this paper.
\begin{theorem}\label{th.uniq}
Let $m$ and $\sigma$ be as in~\eqref{range.exp}. For any $a\in(0,\infty)$, there exist $R(a)\in (0,\infty)$ and a unique positive and decreasing solution $V(\cdot;a)$ to the Cauchy problem~\eqref{CPV} on $(0,R(a))$, and the behavior of $V(\cdot;a)$ as $r\to R(a)$ may be classified as follows.
\begin{enumerate}
  \item There is a unique $a_*\in(0,\infty)$ such that $V(\cdot;a_*)$ has an interface point $R(a_*)$; that is,  $V(R(a_*);a_*) = V'(R(a_*);a_*) = 0$;
  \item For any $a\in(a_*,\infty)$, $V(R(a);a)=0$ and $V'(R(a);a)<0$;
  \item For any $a\in(0,a_*)$, $V(R(a);a)>0$, $V'(R(a);a)=0$ and $V''(R(a);a)>0$. 
\end{enumerate}
\end{theorem}
Let us emphasize once more that the main novelty of Theorem~\ref{th.uniq} is the uniqueness of a non-negative and non-increasing solution to~\eqref{CPV} presenting an interface. The existence of such a solution is already proved by the authors in \cite{IL25}, employing variational techniques.
Theorem~\ref{th.uniq} also provides a complete classification of solutions to~\eqref{CPV} starting from a positive initial value with a derivative featuring possibly a singularity which is weaker than $r^{1-N}$. When $\sigma=0$, a similar uniqueness result is shown in \cite[Section~IV.1.2]{S4} in one space dimension and in \cite{CEF98} when $N\ge 2$.

\medskip

An immediate consequence of Theorem~\ref{th.uniq} and \cite{IL25} is a uniqueness result for~\eqref{eq1} within a specific class of solutions.
\begin{corollary}\label{cor.uni}
Let $m$ and $\sigma$ be as in~\eqref{range.exp}. Then $x\mapsto V(|x|;\blm{a_*})$ is the only non-negative radially symmetric solution to~\eqref{eq1} in $\mathcal{X}$ with non-increasing profile.
\end{corollary}
Indeed, it follows from~\eqref{eq1} and \cite{IL25} that any non-negative radially symmetric solution $w$ to~\eqref{eq1} in $\mathcal{X}$ with non-increasing profile $W$ is a bounded solution to~\eqref{CPV} with compact support; that is, $W=V(\cdot;a_*)$ in view of Theorem~\ref{th.uniq}.

We next report two further applications of Theorem~\ref{th.uniq}, which we have already described in the Introduction. On the one hand, a remarkable application of Theorem~\ref{th.uniq} is the uniqueness (up to scaling) of the extremal function for the CKN inequality \cite[Theorem]{CKN84}
\begin{equation*}
	\mathcal{S}_* N_\sigma^2(w) \le \|\nabla v\|_2^{2a} \|v\|_{(m+1)/m}^{2(1-a)}, \quad v\in \mathcal{X},
\end{equation*}
with $a$ and $\mathcal{S}_*$ defined in~\eqref{ckne} and~\eqref{ckni}, respectively.
\begin{theorem}\label{th.ckn01}
Let $m$ and $\sigma$ be as in~\eqref{range.exp} and consider a non-negative minimizer $w\in \mathcal{X}$, $w\ne 0$, of~\eqref{ckni}; that is,
\begin{equation*}
	w\ge 0 \;\text{ and }\; \mathcal{S}(w) = \mathcal{S}_*. 
\end{equation*} 
Then there are $(\zeta,\eta)\in (0,\infty)^2$ such that $w(x) = \zeta V(\eta |x|;a_*)$, $x\in\mathbb{R}^N$, with $V(\cdot;a_*)$ given by Theorem~\ref{th.uniq}.
\end{theorem}

We next turn to the quasilinear Hardy-H\'enon parabolic equation with source~\eqref{eq.parab} and report the following consequence of Theorem~\ref{th.uniq}.
\begin{theorem}\label{th.parab}
Let $m$ and $\sigma$ be as in~\eqref{range.exp}. Given $T>0$, there is a unique non-negative radially symmetric separate variables solution $u$ in the form~\eqref{SSS} to Eq.~\eqref{eq.parab} with non-increasing profile $U$ such that $u(0)\in\mathcal{X}$. It is given by
\begin{equation*}
	u(t,x):=(T-(m-1)t)^{-1/(m-1)}V^{1/m}(|x|;a_*), \quad (t,x)\in \left[ 0,\frac{T}{m-1} \right)\times\mathbb{R}^N.
\end{equation*}
\end{theorem}
Indeed, we simply observe that, given a non-negative radially symmetric separate variables solution in the form~\eqref{SSS} to Eq.~\eqref{eq.parab} with non-increasing profile $U$, the function $x\mapsto \blm{U^m}(|x|)$ is a non-negative radially symmetric solution to~\eqref{eq1} in $\mathcal{X}$ with non-increasing profile and thus $\blm{U^m}=V(\cdot;a_*)$ by Corollary~\ref{cor.uni}.

\medskip

The proof of Theorem~\ref{th.uniq} is based on studying in parallel the properties of $V(\cdot;a)$ and of the associated variation function
\begin{equation}\label{varfun}
	\varphi(\cdot;a):= \partial_a V(\cdot;a), \quad a\in(0,\infty),
\end{equation}
also considered in papers such as \cite{CEF96, ET98, FQTY08}. The function $\varphi(\cdot;a)$ solves the equation obtained (at least formally) by taking derivatives with respect to $a$ in~\eqref{eq1.rad}; that is,
\begin{subequations}\label{CPphi}
\begin{equation}\label{eq1.phi}
	\varphi''(r)+\frac{N-1}{r}\varphi'(r)+r^{\sigma}\varphi(r)-\frac{1}{m}V^{(1-m)/m}(r;a)\varphi(r)=0, \quad r>0, 
\end{equation}
supplemented with the initial conditions
\begin{equation}\label{icphi}
	\varphi(0)=1, \quad \lim\limits_{r\to0}r^{N-1}\varphi'(r)=0.
\end{equation}
\end{subequations}
Let us recall here that, due to the non-Lipschitz nonlinearity and the singular coefficient $r^{\sigma}$, the well-posedness of both Cauchy problems~\eqref{CPV} and \eqref{CPphi}, as well as the $C^1$-dependence of $V(\cdot;a)$ with respect to $a$ and the rigorous justification of~\eqref{varfun}, do not follow from classical results in the theory of ordinary differential equations. Despite being expected, the proofs of these results are rather involved and technical, adding up to the difficulty of the problem. We believe that the techniques employed in the proofs of these properties can be also useful in the study of other equations involving similar singular terms.

The fundamental step in the proof is the fact that, if $a\in(0,\infty)$ is such that $V(\cdot;a)$ has an interface (in the sense of the first item of Theorem \ref{th.uniq}) at a point $R(a)\in(0,\infty)$, then $\varphi(\cdot;a)$ vanishes exactly once in the interval $(0,R(a))$. Let us remark here that the existence of a unique zero for $\varphi$ is actually a property of (local) monotonicity of the family $V(\cdot;a)$ with respect to $a$ in a neighborhood of a ``good" point: it formally suggests that $V(\cdot;a)$ and $V(\cdot;a')$ should have exactly one intersection point for $a$ and $a'$ close enough. In order to prove this essential property, we adapt to our case the construction and analysis of several auxiliary functions defined and employed in \cite{Tang03, FQTY08}. However, due to the singular variable coefficient, some of the steps in the proof strongly depart from the analogous ones in the above mentioned references.

\medskip

\noindent \textbf{Structure of the paper.} This presentation of the results is followed by Section~\ref{sec.wpae}, in which we gather several properties of the solutions $V(\cdot;a)$ related to their asymptotic expansions as $r\to 0$ and near the interface, if any. 
A comparison principle for Eq.~\eqref{eq1.rad} is stated and proved in Section~\ref{sec.comp}. The core of the proof of Theorem~\ref{th.uniq} is then given in Section~\ref{sec.uniqzero}, where all the needed auxiliary functions are introduced and the existence and uniqueness of a zero for $\varphi(\cdot;a)$ are proved, as discussed in the previous paragraph. The proof of Theorem~\ref{th.uniq} is then completed in Section~\ref{sec.proof}, while Section~\ref{sec.CKN} is devoted to the proof of Theorem~\ref{th.ckn01}. In order to simplify the presentation and not interrupt the main ideas of the proof of Theorem~\ref{th.uniq} with a number of technical results, the (rather lengthy and technical) proofs of the well-posedness of the Cauchy problems~\eqref{CPV} and~\eqref{CPphi}, together with a rigorous justification of~\eqref{varfun}, are postponed to three Appendices.

\section{Well-posedness and asymptotic expansions}\label{sec.wpae}

This section is dedicated to some preparatory results concerning the well-posedness of the Cauchy problems~\eqref{CPV} and~\eqref{CPphi} and the behavior of $V(\cdot;a)$ and of $\varphi(\cdot;a)$, either as $r\to 0$ or as $r$ approaches an interface point $R(a)>0$ of $V(\cdot;a)$, if any, which satisfies $V(R(a);a)=V'(R(a);a)=0$ in this case.

\subsection{Well-posedness of~\eqref{CPV}}\label{sec.wpV}

We begin with the well-posedness of~\eqref{CPV} and, in view of the continuous dependence with respect to $a$ needed in our analysis, we actually study sign-changing solutions to a suitable version of~\eqref{CPV}. Specifically, for $a\in(0,\infty)$, we investigate the well-posedness of the Cauchy problem
\begin{subequations}\label{CPVs} 
\begin{equation}\label{eq1s.rad}
	V''(r)+\frac{N-1}{r}V'(r)+r^{\sigma}V(r)-\big(|V|^{(1-m)/m} V\big)(r)=0, \quad r>0,
\end{equation}
supplemented with the following initial conditions
\begin{equation}\label{icas}
	V(0)=a\in(0,\infty), \quad \lim\limits_{r\to 0}r^{N-1}V'(r)=0.
\end{equation}
\end{subequations}	
In this direction, we report the following result.
\begin{proposition}\label{prop.wpv}
Let $m$ and $\sigma$ be as in~\eqref{range.exp} and $a\in (0,\infty)$. There are $S(a)\in (0,\infty]$ and a unique solution 
\begin{equation*}
	V(\cdot;a)\in C([0,S(a))) \cap C^2((0,S(a)))
\end{equation*}
to~\eqref{CPVs} on $(0,S(a))$ such that
\begin{equation}
	\big(V(r;a),V'(r;a)\big)\ne (0,0), \quad r\in (0,S(a)). \label{nd}
\end{equation}
Moreover,
\begin{itemize}
	\item [$\triangleright$] \textbf{Extension:} if $S(a)<\infty$, then 
	\begin{equation*}
		\lim\limits_{r\to S(a)} V(r;a) = \lim\limits_{r\to S(a)} V'(r;a) = 0
	\end{equation*} 
	and $V(\cdot;a)$ can be extended to a global $C^2$-smooth solution to~\eqref{CPVs} on $(0,\infty)$ by setting $V(r;a):= 0$ for $r\in [S(a),\infty)$;
	\item [$\triangleright$] \textbf{Positivity:} introducing 
	\begin{equation}
		R(a) := \sup\{ r\in (0,S(a))\ :\ V(s;a)>0 \ {\rm and} \ V'(s;a)<0 \ {\rm for} \ s\in(0,r)\}, \label{Ra}
	\end{equation}
	there holds $R(a)\in (0,\infty)$;
	\item [$\triangleright$] \textbf{Continuation:} if $\big(V(R(a);a),V'(R(a);a)\big)\ne (0,0)$, then $R(a)<S(a)$.
\end{itemize}
\end{proposition}

Observe that $(0,R(a))$ is the maximal interval starting from the origin on which $V(\cdot;a)$ is non-negative and decreasing, so that the well-posedness of~\eqref{CPV} within the class of positive and decreasing functions stated in Theorem~\ref{th.uniq} readily follows from Proposition~\ref{prop.wpv}. It is also worth mentioning here that we do not know whether the extension by zero on $(S(a),\infty)$ of $V(\cdot;a)$ when $S(a)<\infty$ is the unique one, a fact which explains why the uniqueness result is restricted to $(0,S(a))$ in that case. 
	
The proof of Proposition~\ref{prop.wpv} is rather technical due to the singularity \blm{of~\eqref{eq1s.rad}} as $r\to 0$ and when $V$ vanishes, and requires several steps. We thus postpone it to Appendix~\ref{app.B}, see Proposition~\ref{prop.CP}, except for the positivity and finiteness of $R(a)$ for any $a>0$\blm{,} which we prove now. 

\begin{lemma}\label{lem.Ra}
For any $a\in (0,\infty)$, $R(a)\in (0,\infty)$. 
\end{lemma}

\begin{proof}
Let $a>0$. Since $V(0;a)=a>0$ and 
\begin{equation}\label{main.id}
	(r^{N-1} V'(r;a))' = -r^{N-1} \left( r^\sigma V(r;a) - \big(|V|^{(1-m)/m}V\big)(r;a) \right) \sim - a r^{\sigma+N-1}, \quad {\rm as} \ r\to 0,
\end{equation}
according to~\eqref{CPVs}, the function $V(\cdot;a)$ is positive and decreasing in a right neighborhood of $r=0$ and thus $R(a)$ is well-defined and positive.

We next assume for contradiction that $R(a)=\infty$. Then
\begin{equation}
	0 < V(r;a) < a \;\;\text{ and }\;\; V'(r;a)<0, \quad r\in (0,\infty). \label{r01}
\end{equation}
In particular, there is $V_\infty\in [0,a)$ such that $V(r;a)\longrightarrow V_\infty$ as $r\to\infty$. \blm{Introducing} the function $T_0\in C^2((0,\infty))$ defined by
\begin{equation*}
	T_0(r) := 1 - r^\sigma V^{(m-1)/m}(r;a), \quad r\in (0,\infty), 
\end{equation*}
the negativity of $\sigma$ and the boundedness of $V(\cdot;a)$ imply that there is $r_*>0$ such that 
\begin{equation}
	\lim_{r\to\infty} T_0(r) = 1 \;\;\text{ and }\;\; T_0(r) > \frac{1}{2}, \quad r\ge r_*. \label{r02}
\end{equation}
Since
\begin{equation}
	\big(r^{N-1} V'(r;a)\big)' = r^{N-1} V^{1/m}(r;a) T_0(r), \quad r\in (0,\infty), \label{r03}
\end{equation}
by~\eqref{eq1.rad}, integrating  the above identity over $(r_*,r)$, $r\ge r_*$, and using~\eqref{r01} give
\begin{equation*}
	- r_*^{N-1} V'(r_*;a) > r^{N-1} V'(r;a) - r_*^{N-1} V'(r_*;a) = \int_{r_*}^r s^{N-1} V^{1/m}(s;a) T_0(s)\, ds \ge 0,
\end{equation*}
from which we deduce that $s \mapsto s^{N-1} V^{1/m}(s;a) T_0(s)$ belongs to $L^1(r_*,\infty)$, and so does $V^{1/m}(\cdot;a) T_0$ due to $r_*>0$. Consequently, the limit of this function as $s\to\infty$, which exists and is equal to $V_\infty^{1/m}$, vanishes and we conclude that
\begin{equation}
	V_\infty = \lim_{r\to\infty} V(r;a) = 0. \label{r04}
\end{equation}

Another consequence of~\eqref{r02} and~\eqref{r03} is that $\big(r^{N-1} V'(r;a)\big)' \ge r^{N-1} V^{1/m}(r;a)/2$ for $r\ge r_*$, hence
\begin{equation}
	-V''(r;a) - \frac{N-1}{r} V'(r;a) + \frac{V^{1/m}(r;a)}{2} < 0, \quad r\in (r_*,\infty). \label{r05}
\end{equation}
At this stage, we adapt an argument from the proof of \cite[Proposition~3.3]{IL25} and introduce the function $\Sigma\in C^2([0,\infty))$ defined by
\begin{equation*}
	\Sigma(r) := (A_*-B_* r^2)_+^{2m/(m-1)}, \quad r\in (0,\infty),
\end{equation*}
with
\begin{equation*}
	A_* := B_* r_*^2 + 2V^{(m-1)/(2m)}(r_*\blm{;a}) \:\;\text{ and }\;\; B_* \left( B_* r_*^2 + 2V^{(m-1)/(2m)}(r_*\blm{;a}) \right) = \frac{(m-1)^2}{16m(m+1)}.
\end{equation*}
Observe that $\Sigma(r_\Sigma)=\Sigma'(r_\Sigma)=\Sigma''(r_\Sigma)=0$ since $2m/(m-1)>2$ and $\Sigma\equiv 0$ on $(r_\Sigma,\infty)$ with $r_\Sigma=\sqrt{A_*/B_*}>r_*$. Also, $\Sigma(r_*)=2V(r_*;a)$. According to \cite[Lemma~3.2]{IL25}, the function $\Sigma$ is a supersolution to~\eqref{r05} and actually satisfies
\begin{equation}
	-\Sigma''(r) - \frac{N-1}{r} \Sigma'(r) + \frac{\Sigma^{1/m}(r)}{2} \ge 0, \quad r\in (0,\infty). \label{r06}
\end{equation}
Consider now $\rho_*$ in $[r_*,r_\Sigma]$ such that 
\begin{equation*}
	W(\rho_*) = \max_{r\in [r_*,r_\Sigma]} W(r) \;\;\text{ with }\;\; W := V(\cdot;a)-\Sigma.
\end{equation*} 
On the one hand, $W'(r_\Sigma)=V'(r_\Sigma;a)<0$ by~\eqref{r01}, so that there is $\delta>0$ such that $W(r)>W(r_\Sigma)=V(r_\Sigma\blm{;a})>0$ for $r\in [r_\Sigma-\delta,r_\Sigma]$. Therefore, $\rho_*<r_\Sigma$ and $W(\rho_*)>0$. On the other hand, $W(r_*)=-V(r_*;a)<0$, from which we deduce that $\rho_*>r_*$. Consequently, $W$ attains a positive maximum at $\rho_*\in (r_*,r_\Sigma)$, which implies that
\begin{equation}
	V(\rho_*;a)>\Sigma(\rho_*), \quad V'(\rho_*;a)=\Sigma'(\rho_*) \;\;\text{ and }\;\; V''(\rho_*;a) \le \Sigma''(\rho_*). \label{r07}
\end{equation}
It then follows from~\eqref{r05}, \eqref{r06} and~\eqref{r07} that
\begin{align*}
	0 & > - V''(\rho_*;a) - \frac{N-1}{\rho_*} V'(\rho_*;a) + \frac{V^{1/m}(\rho_*;a)}{2} + \Sigma''(\rho_*) + \frac{N-1}{\rho_*} \Sigma'(\rho_*) - \frac{\Sigma^{1/m}(\rho_*)}{2} \\
	& = \Sigma''(\rho_*) - V''(\rho_*;a) + \frac{\big(V^{1/m}(\rho_*;a)-\Sigma^{1/m}(\rho_*)\big)}{2} \ge 0,
\end{align*}
and a contradiction. Therefore, $R(a)$ is finite and the proof is complete.
\end{proof}

\begin{proof}[Proof of Proposition~\ref{prop.wpv}]
Let $a>0$. The existence of $S(a)$ and the well-posedness of~\eqref{CPVs} on $(0,S(a))$ follows from Proposition~\ref{prop.CP}, while the positivity and finiteness of $R(a)$ are provided by Lemma~\ref{lem.Ra}. Furthermore, the definition of $R(a)$ guarantees that $R(a)\le S(a)$. \blm{Now,} either $S(a)=\infty$ and Lemma~\ref{lem.Ra} readily implies that $R(a)<S(a)$, or $S(a)<\infty$ and, if
\blm{\begin{equation*}
\big(V(R(a);a),V'(R(a);a)\big)\ne (0,0)=\big(V(S(a);a),V'(S(a);a)\big),
\end{equation*}
then} necessarily $R(a)\ne S(a)$ and thus $R(a)<S(a)$.
\end{proof}

\subsection{Well-posedness of~\eqref{CPphi}}\label{sec.wpp}

We next turn to the Cauchy problem~\eqref{CPphi}, which is also well-posed on $(0,R(a))$ for any $a\in (0,\infty)$, according to Lemma~\ref{lem.exist} with $R=R(a)$, $U=V^{(m-1)/m}(\cdot;a)/m$ and $h_0=1$. We refer to Proposition~\ref{prop.wpvphib} below for a complete proof but report already here the well-posedness result for clarity.
\begin{proposition}\label{prop.wpvphi}
Let $m$ and $\sigma$ be as in~\eqref{range.exp} and $a\in (0,\infty)$. Then, there is a unique solution $\varphi(\cdot;a)\in C([0,R(a))) \cap C^2((0,R(a)))$ to~\eqref{CPphi} and
\begin{equation}
	\varphi'(r;a) = - \frac{r^{\sigma+1}}{\sigma+N} + o(r^{\sigma+1}), \quad {\rm as} \ r\to 0. \label{phip0}
\end{equation}
\end{proposition}

\subsection{Asymptotic expansion as $r\to0$}\label{sec.ae0}

The next result deals with the expansion of $V(\cdot;a)$ near zero, and is a simple consequence of the analysis performed in \cite[Theorem~1.2]{IL25}.
\begin{lemma}\label{lem.exp0}
Let $m$ and $\sigma$ be as in~\eqref{range.exp}. Then, for any $a\in(0,\infty)$ we have, as $r\to 0$,
\begin{equation}\label{exp0}
	V'(r;a)=-\frac{a}{\sigma+N}r^{\sigma+1}+o(r^{\sigma+1}), \quad V(r;a)=a-\frac{a}{(\sigma+N)(\sigma+2)}r^{\sigma+2}+o(r^{\sigma+2}).
\end{equation}
\end{lemma}
\begin{proof}
This result is proved in \cite[Section~3.3]{IL25}. We give below a sketch of the proof for the sake of completeness. We start from the identity~\eqref{main.id} and recall that $V(0;a)=a$. We thus deduce that
\begin{equation*}
	(r^{N-1}V'(r))' = -a r^{N+\sigma-1}+o(r^{N+\sigma-1}), \quad {\rm as} \ r\to 0,
\end{equation*}
which, along with the initial condition \eqref{ica} and an application of L'H\^{o}pital's rule\blm{,} gives
\begin{equation*}
	\lim\limits_{r\to 0}\frac{(\sigma+N)r^{N-1}V'(r)}{r^{N+\sigma}} = -a\blm{.}
\end{equation*}
\blm{The latter limit} is equivalent to the expansion of $V'(\cdot;a)$ in~\eqref{exp0}. The expansion of $V(\cdot;a)$ follows then by integration.
 \end{proof}
 
 \subsection{Asymptotic expansion near the interface}\label{sec.aei}
 
The asymptotic expansion near the interface point (when such an interface point exists) is a bit more involved and is established in the next result.
\begin{lemma}\label{lem.expint}
Assume that $a>0$ is such that $V(R(a);a)=V'(R(a);a)=0$. Then the following interface behavior holds true:
\begin{equation}\label{expint}
	V(r;a)\sim C_I(R(a)-r)^{2m/(m-1)}, \quad C_I:=\left[\frac{m-1}{\sqrt{2m(m+1)}}\right]^{2m/(m-1)}, \quad {\rm as} \ r\to R(a).
\end{equation}
\end{lemma}
\begin{proof}
We employ an idea inspired by \cite{BPT86}. Since $a>0$ is fixed, we drop for simplicity the dependence on $a$ from the notation. Introducing $H(r):=r^{N-1}V'(r)$ for $r\in [0,R(a)]$, we infer from~\eqref{main.id} that $H$ satisfies
\begin{equation}\label{expint1}
	H'(r)=r^{N-1}V^{1/m}(r)\left(1-r^{\sigma}V^{(m-1)/m}(r)\right), \quad r\in (0,R(a)).
\end{equation}
Multiplying both sides of the identity~\eqref{expint1} by $2H(r)$, $r\in (0,R(a))$, we find
\begin{equation*}
\begin{split}
	(H^2)'(r)&=2r^{2(N-1)}V'(r)V^{1/m}(r)\left(1-r^{\sigma}V^{(m-1)/m}(r)\right)\\
	&=\frac{2m}{m+1}r^{2(N-1)}\big(V^{(m+1)/m}\big)'(r)\left(1-r^{\sigma}V^{(m-1)/m}(r)\right),
\end{split}
\end{equation*}
or, equivalently,
\begin{equation}\label{expint2}
	\frac{\big(H^2\big)'(r)}{\big(V^{(m+1)/m}\big)'(r)}=\frac{2m}{m+1}r^{2(N-1)}\left(1-r^{\sigma}V^{(m-1)/m}(r)\right).
\end{equation}
Since $m>1$, $R(a)>0$ and $V(R(a))=0$, we readily conclude from~\eqref{expint2} that
\begin{equation*}
	\lim\limits_{r\to R(a)}\frac{\big(H^2\big)'(r)}{\big(V^{(m+1)/m}\big)'(r)}=\frac{2m}{m+1}R(a)^{2(N-1)}.
\end{equation*}
Recalling that $H(R(a))=V(R(a))=0$, we obtain from L'H\^{o}pital's rule that
\begin{equation*}
	\lim\limits_{r\to R(a)}\frac{H^2(r)}{V^{(m+1)/m}(r)}=\frac{2m}{m+1}R(a)^{2(N-1)};
\end{equation*}
that is, in view of the negativity of $H$,
\begin{equation}\label{expint3}
	\lim\limits_{r\to R(a)}\frac{-H(r)}{V^{(m+1)/2m}(r)}=\left(\frac{2m}{m+1}\right)^{1/2}R(a)^{N-1}.
\end{equation}
Since
\begin{equation*}
	\frac{-H(r)}{V^{(m+1)/2m}(r)}=-r^{N-1}V'(r)V^{-(m+1)/2m}(r)=-\frac{2m}{m-1}r^{N-1}\big(V^{(m-1)/2m}\big)'(r),
\end{equation*}
we deduce from~\eqref{expint3}, the positivity of $R(a)$ and obvious simplifications that
\begin{equation}\label{expint4}
	\lim\limits_{r\to R(a)}\big(V^{(m-1)/2m}\big)'(r)=-\frac{m-1}{\sqrt{2m(m+1)}}= - C_I^{(m-1)/2m}.
\end{equation}
The interface behavior~\eqref{expint} is then obtained by integration from \eqref{expint4}, completing the proof.
\end{proof}

\section{A comparison principle}\label{sec.comp}

In this section we prove a comparison principle for~\eqref{eq1.rad} which is analogous to \cite[Lemma~4.4]{CEF96} and is in fact a variation of the elliptic comparison principle introduced in \cite{BNV}. However, since~\eqref{eq1.rad} features a variable coefficient, some changes in the statement and proof are needed.
\begin{proposition}\label{prop.comp}
Let $0<r_1<r_2$ and let $g$ and $h$ be two solutions in $C^2([r_1,r_2])$ to
\begin{equation}\label{eq1.rads}
	f''(r)+\frac{N-1}{r}f'(r)+r^{\sigma}f(r)- \big(|f|^{(1-m)/m}f\big)(r)=0, \quad r\in [r_1,r_2].
\end{equation} 
Assume that there is $A>0$ such that
\begin{equation}\label{cond.comp1}
	0\leq g(r)\leq A, \quad r\in(r_1,r_2), \quad h(r_1)<g(r_1), \quad h'(r_1)<g'(r_1),
\end{equation}
and
\begin{equation}\label{cond.comp2}
	z\mapsto r^{\sigma}z- z^{1/m} \ {\rm is \ decreasing \ on} \ (0,A), \quad {\rm for \ all} \ r\in(r_1,r_2).
\end{equation}
Then

\noindent (a) if $h(r)\geq0$ for $r\in[r_1,r_2]$, then $h(r)<g(r)$ for $r\in[r_1,r_2]$ and $h'(r_2)<g'(r_2)$;

\noindent (b) if $h(r_1)>0$ and $(gg')(r_2)=0$, then there is $r_0\in(r_1,r_2)$ such that $h>0$ on $[r_1,r_0)$, $h(r_0)=0$ and $h'(r_0)<0$.
\end{proposition}
\begin{proof}
\noindent (a) Assume for contradiction that there is $\bar{r}\in(r_1,r_2]$ such that $h(r)<g(r)$ for $r\in(r_1,\bar{r})$ and $g(\bar{r})=h(\bar{r})$. By subtracting~\eqref{eq1.rad} satisfied by $g$ and $h$, we readily deduce that, for $r\in(r_1,r_2)$,
\begin{equation*}
\begin{split}
	\frac{d}{dr}\big[r^{N-1}(g-h)'(r)\big]&=r^{N-1}\left[(g-h)''(r)+\frac{N-1}{r}(g-h)'(r)\right]\\
	&=r^{N-1}\left[-r^{\sigma}(g-h)(r)+(g^{1/m}-h^{1/m})(r)\right].
\end{split}
\end{equation*}
We multiply by $g-h$ and integrate over $(r_1,\bar{r})$ to get
\begin{equation}\label{comp1}
\begin{split}
	& \int_{r_1}^{\bar{r}}\big[r^{N-1}(g-h)'(r)\big]'(g-h)(r)\,dr \\
	& \hspace{1cm} =\int_{r_1}^{\bar{r}} r^{N-1}(g-h)(r) \left[-r^{\sigma}(g-h)(r)+(g^{1/m}-h^{1/m})(r)\right]\,dr.
\end{split}
\end{equation}
An integration by parts in the left hand side of~\eqref{comp1}, together with the definition of $\bar{r}$, gives
\begin{equation}\label{comp2}
\begin{split}
	-r_1^{N-1}(g'-h')(r_1)&(g-h)(r_1) -\int_{r_1}^{\bar{r}}r^{N-1}((g-h)')^2(r)\,dr\\
&= - \int_{r_1}^{\bar{r}} r^{N-1}\left[r^{\sigma}(g-h)(r)-(g^{1/m}-h^{1/m})(r) \right] (g-h)(r)\,dr.
\end{split}
\end{equation}
On the one hand, we infer from~\eqref{cond.comp1} that the left hand side of~\eqref{comp2} is negative. On the other hand, we deduce from~\eqref{cond.comp1} and the definition of $\bar{r}$ that $0\leq h(r)\leq g(r)\leq A$ for $r\in[r_1,\bar{r}]$, and the monotonicity assumption ~\eqref{cond.comp2} implies that the right hand side of~\eqref{comp2} is non-negative. We have thus reached a contradiction, which stems from the assumption of the existence of $\bar{r}$. We have therefore proved that $h(r)<g(r)$ for $r\in[r_1,r_2]$. Assume next for contradiction that $g'(r_2)\le h'(r_2)$. We then proceed as in the proof of~\eqref{comp2} to obtain
\begin{equation}\label{comp3}
\begin{split}
	r_2^{N-1} (g'-h')(r_2)(g-h)(r_2) & - r_1^{N-1} (g'-h')(r_1) (g-h)(r_1) -\int_{r_1}^{r_2} r^{N-1} ((g-h)')^2(r)\,dr\\
	&= - \int_{r_1}^{r_2} r^{N-1}\left[r^{\sigma}(g-h)(r)-(g^{1/m}-h^{1/m})(r) \right] (g-h)(r)\,dr.
\end{split}
\end{equation}
Since $(g-h)(r_2)>0$ and $(g'-h')(r_2)\le 0$, the left hand side of~\eqref{comp3} is negative and we are led to an analogous contradiction, completing the proof of part~(a).

\medskip

\noindent (b) Assume for contradiction that $h>0$ on $(r_1,r_2)$. Then, $h\ge 0$ on $[r_1,r_2]$, so that, according to part~(a), we deduce that $h(r)<g(r)$ for $r\in[r_1,r_2]$. In particular, $h(r_2)<g(r_2)=0$ and we reach a contradiction. Therefore, $h$ vanishes at least once in $(r_1,r_2)$ and the positivity of $h(r_1)$ implies that there is $r_0\in(r_1,r_2)$ such that $h(r)>0$ for $r\in [r_1,r_0)$ and $h(r_0)=0$. In particular, $h'(r_0)\le 0$. Now, if $h'(r_0)=0=h(r_0)$, then we have also $h''(r_0)=0$ by~\eqref{eq1.rads}, so that the function $\bar{h}$ defined by $\bar{h}(r)=h(r)$ for $r\in [r_1,r_0]$ and $\bar{h}(r)=0$ for $r\in [r_0,r_2]$ belongs to $C^2([r_1,r_2])$ and also solves~\eqref{eq1.rads} on $(r_1,r_2)$. Since $\bar{h}\ge 0$ on $[r_1,r_2]$ and satisfies~\eqref{cond.comp1}, we infer from part~(a) that $0=\bar{h}(r_2)<g(r_2)$ and $0=\bar{h}'(r_2)<g'(r_2)$. Therefore, $(gg')(r_2)>0$ and a contradiction. We have thus established that $h'(r_0)<0$ and thereby completed the proof of part~(b).
\end{proof}

\section{Uniqueness of a zero of $\varphi(\cdot;a)$}\label{sec.uniqzero}

The aim of this section is to establish the fact that, if $a\in (0,\infty)$ is such that $V(R(a);a)=V'(R(a);a)=0$, then $\varphi(\cdot;a)$ has a unique zero in the interval $(0,R(a))$. In order to prove this result, we follow a general scheme stemming from \cite{Tang03}, also employed in \cite{FQTY08}. Let us stress that, in these references, the equation under study \blm{is} autonomous and \blm{has} no singular \blm{coefficient}, and the adaptation to cope with the singular coefficient $r^{\sigma}$ introduces some technical complications. We start by introducing a number of auxiliary functions and study their properties. Throughout this section, in order to simplify the notation \blm{and} whenever there is no risk of confusion, we drop the dependence on $a$ and write $V$ and $\varphi$ instead of $V(\cdot;a)$ and $\varphi(\cdot;a)$.
\begin{lemma}\label{lem.energy}
Let $a\in (0,\infty)$. Define
\begin{equation}\label{energy}
	E(r;a):=\frac{1}{2}(V')^2(r;a)+\frac{1}{2}r^{\sigma}V^2(r;a)-\frac{m}{m+1}V^{(m+1)/m}(r;a), \quad r\in (0,R(a)].
\end{equation}
Then $E(\cdot;a)$ is decreasing on $(0,R(a)]$ with $E(r;a)\ge - m a^{(m+1)/m}/(m+1)$ for $r\in (0,R(a)]$ and
\begin{equation*}
	E(r;a) \sim \frac{a^2}{2} r^\sigma,  \quad {\rm as} \ r\to 0.
\end{equation*}.
\end{lemma}
\begin{proof}
Let $r\in (0,R(a)\blm{)}$. A direct calculation, along with the negativity of $\sigma$, gives
\begin{equation*}
\begin{split}
	E'(r)&=V'(r)V''(r)+r^{\sigma}V(r)V'(r)-V^{1/m}(r)V'(r)+\frac{\sigma}{2}r^{\sigma-1}V^2(r)\\
	&=-\frac{N-1}{r}(V')^2(r)+\frac{\sigma}{2}r^{\sigma-1}V^2(r)<0,
\end{split}
\end{equation*}
as claimed. The lower bound readily follows from the properties of $V$ due to the definition of $R(a)$, while the behavior as $r\to 0$ is an immediate consequence of~\eqref{exp0} and the positivity of $\sigma+2$.
\end{proof}
As in \cite{Tang03}, we next introduce two auxiliary functions that will appear in the proof of the uniqueness of the zero of $\varphi$. Set, for $r\in(0,R(a))$,
\begin{equation}\label{fPQ}
	P(r):=2r^NE(r)+(N-2)r^{N-1}V(r)V'(r), \quad Q(r):=P(r)+\frac{m-1}{m+1}r^NV^{(m+1)/m}(r).
\end{equation}
\begin{lemma}\label{lem.PQ}
Let $a\in(0,\infty)$ be such that $V(R(a))=0$. Then $P(r)>0$ and $Q(r)>0$ for any $r\in(0,R(a))$.
\end{lemma}
\begin{proof}
Let $r\in(0,R(a))$. We note that
\begin{equation*}
	P(r)=r^N(V')^2(r)+r^{\sigma+N}V^2(r)-\frac{2m}{m+1}r^{N}V^{(m+1)/m}(r)+(N-2)r^{N-1}V(r)V'(r).
\end{equation*}
Making use of the previous expression and~\eqref{eq1.rad}, we compute
\begin{equation*}
\begin{split}
	P'(r)&=Nr^{N-1}(V')^2(r)+2r^NV'(r)\left(-\frac{N-1}{r}V'(r)-r^{\sigma}V(r)+V^{1/m}(r)\right)\\
	&+(\sigma+N)r^{\sigma+N-1}V^2(r)+2r^{\sigma+N}V(r)V'(r)-\frac{2mN}{m+1}r^{N-1}V^{(m+1)/m}(r)-2r^{N}V^{1/m}(r)V'(r)\\
	&+(N-1)(N-2)r^{N-2}V(r)V'(r)+(N-2)r^{N-1}(V')^2(r)\\
	&+(N-2)r^{N-1}V(r)\left(-\frac{N-1}{r}V'(r)-r^{\sigma}V(r)+V^{1/m}(r)\right),
\end{split}
\end{equation*}
and, after simplifying similar terms, we are left with
\begin{equation*}
	P'(r)=(\sigma+2)r^{\sigma+N-1}V^2(r)+\frac{(N-2)(m+1)-2mN}{m+1}r^{N-1}V^{(m+1)/m}(r),
\end{equation*}
which can be written in the equivalent, but more convenient, form
\begin{equation}\label{deriv.P}
	P'(r)=(\sigma+2)r^{\sigma+N-1}V^{(m+1)/m}(r)\left(V^{(m-1)/m}(r)-\frac{N(m-1)+2(m+1)}{(m+1)(\sigma+2)}r^{-\sigma}\right),
\end{equation}
for $r\in(0,R(a))$. The definition~\eqref{Ra} of $R(a)$, the negativity of $\sigma$ and the fact that $m>1$ imply that the mapping
\begin{equation*}
	P_1(r):=V^{(m-1)/m}(r)-\frac{N(m-1)+2(m+1)}{(m+1)(\sigma+2)}r^{-\sigma}, \quad r\in (0,R(a)),
\end{equation*}
is decreasing on $(0,R(a))$, while the assumption that $V(R(a))=0$ gives
\begin{equation*}
	\lim\limits_{r\to 0} P_1(r) = a^{(m-1)/m}>0, \quad \lim\limits_{r\to R(a)} P_1(r) = -\frac{N(m-1)+2(m+1)}{(m+1)(\sigma+2)} R(a)^{-\sigma} < 0.
\end{equation*}
Since $\sigma+2>0$, it follows from~\eqref{deriv.P} that there exists a unique point $r_1(a)\in(0,R(a))$ such that $P'(r_1(a))=0$ and
\begin{equation*}
	P'(r)>0, \quad {\rm for} \ r\in(0,r_1(a)), \quad P'(r)<0, \quad {\rm for} \ r\in(r_1(a),R(a)).
\end{equation*}
On the one hand, recalling the expansion of $V'$ as $r\to 0$ given in Lemma~\ref{lem.exp0}, we have
\begin{equation*}
	r^{N-1}(VV')(r)\sim-\frac{a^2}{\sigma+N}r^{\sigma+N}, \quad r^N(V')^2(r)\sim\frac{a^2}{(\sigma+N)^2}r^{N+2\sigma+2}, \quad {\rm as} \ r\to0,
\end{equation*}
and the positivity of $N+\sigma>0$, $N+2\sigma+2=(N+\sigma)+(\sigma+2)>0$, along with~\eqref{ica}, ensures that $P(r)\to0$ as $r\to0$. On the other hand, the assumption $V(R(a))=0$ implies that
\begin{equation*}
	\lim\limits_{r\to R(a)}P(r)=R(a)^N(V'(R(a)))^2\geq0.
\end{equation*}
The uniqueness of the critical point $r_1(a)\in(0,R(a))$ then yields $P>0$ on $(0,R(a))$ and also $Q(r)>P(r)>0$ for $r\in(0,R(a))$, concluding the proof.
\end{proof}
The next step is to define two new auxiliary functions, the first being the Wronskian between $V$ and $\varphi$ for the same value of $a$:
\begin{subequations}\label{xi.zeta}
\begin{equation}\label{fxi}
	\xi(r):=r^{N-1} (V'\varphi-V\varphi')(r), \quad r\in(0,R(a)),
\end{equation}
\begin{equation}\label{fzeta}
	\zeta(r):=r^{N}\left[ (V'\varphi')(r)+(r^{\sigma}V(r)-V^{1/m}(r))\varphi(r)\right]+(N-2)r^{N-1} (V'\varphi)(r), \quad r\in(0,R(a)).
\end{equation}
\end{subequations}
We compute the first derivatives of these functions in the next result.
\begin{lemma}\label{lem.xizeta}
We have
\begin{equation}\label{deriv.xi}
	\xi'(r)=\frac{m-1}{m}r^{N-1}V^{1/m}(r)\varphi(r), \quad r\in(0,R(a)),
\end{equation}
and
\begin{equation}\label{deriv.zeta}
	\zeta'(r)=(\sigma+2)r^{N-1+\sigma}V^{1/m}(r)\varphi(r)\left[V^{(m-1)/m}(r)-\frac{2}{\sigma+2}r^{-\sigma}\right], \quad r\in(0,R(a)).
\end{equation}
Moreover,
\begin{equation}\label{exp.xizeta}
	\lim\limits_{r\to 0}\xi(r)=\lim\limits_{r\to 0}\zeta(r)=0.
\end{equation}
\end{lemma}
\begin{proof}
We proceed by direct calculations, employing the equations~\eqref{eq1.rad} solved by $V$ and~\eqref{eq1.phi} solved by $\varphi$, noticing that some great simplifications take place. Indeed, we have
\begin{equation*}
\begin{split}
\xi'(r)&=(N-1)r^{N-2}(V'\varphi-V\varphi')(r)+r^{N-1}(V''\varphi-V\varphi'')(r) \\
     &=(N-1)r^{N-2}(V'\varphi-V\varphi')(r)+r^{N-1}\varphi(r)\left(-\frac{N-1}{r}V'(r)-r^{\sigma}V(r)+V^{1/m}(r)\right)\\
     &-r^{N-1}V(r)\left(-\frac{N-1}{r}\varphi'(r)-r^{\sigma}\varphi(r)+\frac{1}{m}V^{(1-m)/m}(r)\varphi(r)\right)\\
     &=r^{N-1}V^{1/m}(r)\varphi(r)-\frac{1}{m}r^{N-1}V^{1/m}(r)\varphi(r),
\end{split}
\end{equation*}
which implies~\eqref{deriv.xi}. A similar (but longer) calculation gives
\begin{equation*}
\begin{split}
\zeta'(r)& = N r^{N-1} V'(r) \varphi'(r) + r^N \varphi'(r) \left( -\frac{N-1}{r} V'(r) - r^{\sigma} V(r) + V^{1/m}(r) \right) \\
& + r^{N-1} V'(r) \left( -\frac{N-1}{r} \varphi'(r) - r^{\sigma} \varphi(r) + \frac{1}{m} V^{(1-m)/m}(r) \varphi(r) \right) + N r^{N-1} \big(r^{\sigma} V(r) - V^{1/m}(r)\big) \varphi(r) \\
& + r^{N} \left( \sigma r^{\sigma-1} V(r) + r^{\sigma} V'(r) - \frac{1}{m} V^{(1-m)/m}(r) V'(r)\right) \varphi(r)+ r^{N} \big(r^{\sigma} V(r) - V^{1/m}(r)\big)\varphi'(r)\\
& + (N-2)(N-1) r^{N-2} V'(r) \varphi(r) + (N-2) r^{N-1} V'(r) \varphi'(r) \\
& + (N-2) r^{N-1} \varphi(r) \left( -\frac{N-1}{r} V'(r) - r^{\sigma} V(r) + V^{1/m}(r) \right).
\end{split}
\end{equation*}
After simplifying similar terms in the long expression of $\zeta'(r)$, we are only left with
\begin{equation*}
\zeta'(r)=(\sigma+2)r^{\sigma+N-1}V(r)\varphi(r)-2r^{N-1}V^{1/m}(r)\varphi(r),
\end{equation*}
which is equivalent to~\eqref{deriv.zeta}. We are left to prove~\eqref{exp.xizeta}, which follows easily from Lemma~\ref{lem.exp0}. Indeed, taking into account the expansion~\eqref{exp0} of $V'$, the expansion~\eqref{phip0} satisfied by $\varphi'$ and the definition~\eqref{fxi}, we find
\begin{equation*}
	\xi(r)=r^{N-1}\left(-\frac{a}{\sigma+N}r^{\sigma+1}+o(r^{\sigma+1})+\frac{a}{\sigma+N}r^{\sigma+1}+o(r^{\sigma+1})\right)=o(r^{\sigma+N})
\end{equation*}
and the positivity of $\sigma+N$ implies that $\xi(r)\to0$ as $r\to0$. Using the same expansions as above and the definition~\eqref{fzeta} of the function $\zeta$ gives
\begin{equation*}
\begin{split}
	\zeta(r)&=r^{N}\left[\frac{a}{(\sigma+N)^2}r^{2(\sigma+1)}+o(r^{2(\sigma+1)})\right]+ar^{N+\sigma}+o(r^{N+\sigma}) \\
    &-a^{1/m}r^N+o(r^N)-\frac{(N-2)a}{\sigma+N}r^{\sigma+N}+o(r^{\sigma+N}) \\
    &=\frac{a(\sigma+2)}{\sigma+N}r^{\sigma+N}+o(r^{\sigma+N}),
\end{split}
\end{equation*}
and once more the positivity of $\sigma+N$ ensures that $\zeta(r)\to0$ as $r\to 0$, concluding the proof of~\eqref{exp.xizeta}.
\end{proof}
Based on the definitions~\eqref{xi.zeta} of the functions $\xi$ and $\zeta$ and following \cite{Tang03, FQTY08}, we define our last auxiliary function, which is closely related to the Wronskian of $\xi$ and $\zeta$: for $r\in(0,R(a))$, \blm{set}
\begin{equation}\label{fT}
\begin{split}
	T(r)& := \frac{\zeta'(r)\xi(r)-\zeta(r)\xi'(r)}{\xi'(r)}\\
	&=\frac{m(\sigma+2)}{m-1}r^{\sigma}\left[V^{(m-1)/m}(r)-\frac{2}{\sigma+2}r^{-\sigma}\right]\xi(r)-\zeta(r),
\end{split}
\end{equation}
\blm{the second equality being ensured by} Lemma~\ref{lem.xizeta}. We gather basic properties of the function $T$ in the next result.
\begin{lemma}\label{lem.T}
For $r\in (0,R(a))$,
\begin{equation}\label{deriv.T}
	T'(r)=\frac{m(\sigma+2)}{m-1}r^{\sigma}V^{-1/m}(r)\left[\frac{\sigma}{r}V(r)+\frac{m-1}{m}V'(r)\right]\xi(r)
\end{equation}
and
\begin{equation}\label{exp.T}
	\lim\limits_{r\to0}T(r)=0.
\end{equation}
\end{lemma}
\begin{proof}
We first prove~\eqref{exp.T}. To this end, \blm{taking into account the initial conditions $V(0)=a>0$ and $\varphi(0)=1$, we observe from~\eqref{deriv.xi} that}
\begin{equation*}
	\xi'(r)=\frac{m-1}{m}a^{1/m}r^{N-1}+o(r^{N-1}), \quad {\rm as} \ r\to 0,
\end{equation*}
so that
\begin{equation*}
	\xi(r)=\frac{m-1}{mN}a^{1/m}r^N+o(r^N), \quad {\rm as} \ r\to 0,
\end{equation*}
and thus
\begin{equation}\label{interm17}
	\lim\limits_{r\to0}r^{\sigma}\xi(r)=0.
\end{equation}
Since
\begin{equation*}
	T(r)=\frac{m(\sigma+2)}{m-1}V^{(m-1)/m}(r)r^{\sigma}\xi(r)-\frac{2m}{m-1}\xi(r)-\zeta(r),
\end{equation*}
the property~\eqref{exp.T} follows from~\eqref{exp.xizeta} and~\eqref{interm17}. We next turn our attention to the expression of $T'$. Let $r\in (0,R(a))$. Since
\begin{equation*}
	T(r)=\frac{\zeta'(r)}{\xi'(r)}\xi(r)-\zeta(r),
\end{equation*}
we deduce that
\begin{equation}
\begin{split}
	T'(r)&= \left( \frac{\zeta'}{\xi'} \right)'(r) \xi(r)\\
	&=\frac{m(\sigma+2)}{m-1}\frac{d}{dr}\left\{r^{\sigma}\left[V^{(m-1)/m}(r)-\frac{2}{\sigma+2}r^{-\sigma}\right]\right\}\xi(r)\\
	&=\frac{m(\sigma+2)}{m-1}\frac{d}{dr}\left[r^{\sigma}V^{(m-1)/m}(r)-\frac{2}{\sigma+2}\right]\xi(r),
\end{split}
\end{equation}
which readily leads to \eqref{deriv.T}.
\end{proof}
We are now in a position to proceed to the proof of the main result of this section; that is, the existence and uniqueness of a zero of the function $\varphi$ for any $a\in(0,\infty)$ such that $V(R(a)) =V'(R(a))=0$. This will be the outcome of the following two propositions.
\begin{proposition}\label{prop.existzero}
Let $a>0$ such that $V(R(a))=V'(R(a))=0$. Then $\varphi$ vanishes at least once in $(0,R(a))$.
\end{proposition}
\begin{proof}
We introduce the Wronskian $W:=V\varphi'-\varphi V'$ and \blm{infer from~\eqref{deriv.xi} that}
\begin{equation}\label{interm8}
	\frac{d}{dr}(r^{N-1}W(r)) = \blm{-\xi'(r)} = -\frac{m-1}{m}r^{N-1}V^{1/m}(r)\varphi(r), \quad r\in (0,R(a)).
\end{equation}
The remaining part of the proof is divided in several steps.

\medskip

\noindent \textbf{Step~1.} Assume now for contradiction that $\varphi>0$ on $(0,R(a))$. Under this assumption, we also claim that
\begin{equation}\label{cond.zero}
	\liminf\limits_{r\to R(a)} V(r)\varphi'(r)\geq0.
\end{equation}
Let us assume, in this first step, that~\eqref{cond.zero} is fulfilled. Owing to the assumed positivity of $\varphi$, the right-hand side of~\eqref{interm8} is negative on $(0,R(a))$ and thus the mapping $r\mapsto r^{N-1}W(r)$ is decreasing on $(0,R(a))$. The initial conditions~\eqref{icphi} (satisfied by $\varphi$) and~\eqref{ica} (satisfied by $V$) ensure that
\begin{equation*}
	R(a)^{N-1} \limsup_{r\to R(a)} W(r) \le r^{N-1}W(r) <\lim\limits_{r\to0} r^{N-1}W(r)=0, \quad r\in(0,R(a)).
\end{equation*}
But the definition~\eqref{Ra} of $R(a)$, together with the positivity assumption on $\varphi$, implies that
\begin{equation*}
	V'(r)\varphi(r)<0, \quad r\in(0,R(a)),
\end{equation*}
and thus, in view of~\eqref{cond.zero},
\begin{equation*}
	\limsup_{r\to R(a)} W(r) = \limsup_{r\to R(a)} \big[V(r) \varphi'(r) - V'(r) \varphi(r)\big] \ge \limsup_{r\to R(a)} \blm{\big[V(r) \varphi'(r) \big]} \ge 0. 
\end{equation*}
We thus reach a contradiction, hence $\varphi$ vanishes at least once in $(0,R(a))$.

\medskip

\noindent \textbf{Step 2.} We are left to prove the claim~\eqref{cond.zero}, under the assumption that $\varphi>0$ in $(0,R(a))$, for the previous argument of contradiction and thus the proof to be complete. Since $\sigma<0$ and $V$ is decreasing on $(0,R(a))$, the function
\begin{equation*}
	r\mapsto 1-mr^{\sigma}V^{(m-1)/m}(r), \quad r\in(0,R(a)),
\end{equation*}
is increasing, with 
\begin{equation*}
	\lim_{r\to 0} \big(1-mr^{\sigma}V^{(m-1)/m}(r)\big) = - \infty \;\;\text{ and }\;\; \lim_{r\to R(a)} \big(1-mr^{\sigma}V^{(m-1)/m}(r)\big) = 1.
\end{equation*}
There is thus a unique $\varrho_0\in(0,R(a))$ such that
\begin{equation}
	1-mr^{\sigma}V^{(m-1)/m}(r)\begin{cases}
                             <0, & \mbox{if } r\in(0,\varrho_0), \\
                             >0, & \mbox{if } r\in(\varrho_0,R(a)).
                           \end{cases}\label{ds}
\end{equation}
Owing to the positivity of $V$ and the assumed positivity of $\varphi$ on $(0,R(a))$, we deduce from the previous property that
\begin{equation}\label{interm18}
	(r^{N-1}\varphi'(r))'=\frac{r^{N-1}}{m}V^{(1-m)/m}(r)\varphi(r)\left[1-mr^{\sigma}V^{(m-1)/m}(r)\right]\begin{cases}
                             <0, & \mbox{if } r\in(0,\varrho_0), \\
                             >0, & \mbox{if } r\in(\varrho_0,R(a)),
                           \end{cases}
\end{equation}
which in particular implies that $r\mapsto r^{N-1}\varphi'(r)$ is a decreasing function on $(0,\varrho_0)$. This monotonicity property, along with~\eqref{icphi}, further entails that $\varphi'(r)<0$ for $r\in(0,\varrho_0]$. We split now the analysis into two cases, according to whether $\varphi'$ vanishes or not in the complementary interval $(\varrho_0,R(a))$.

\medskip

\noindent \textbf{Step 3. Case 1:} there is $\varrho_1\in(\varrho_0,R(a))$ such that
\begin{equation*}
	\varphi'(r)<0, \quad r\in(\varrho_0,\varrho_1), \quad \varphi'(\varrho_1)=0.
\end{equation*}
In this case, integrating \eqref{interm18} over $(\varrho_1,r)$ for $r\in(\varrho_1,R(a))$ leads to
\begin{equation*}
	r^{N-1}\varphi'(r)=\int_{\varrho_1}^r\frac{s^{N-1}}{m}V^{(1-m)/m}(s)\varphi(s)\left[1-ms^{\sigma}V^{(m-1)/m}(s)\right]\,ds>0,
\end{equation*}
hence $\varphi'(r)>0$ in $(\varrho_1,R(a))$. The non-negativity of $V$ gives then the claim \eqref{cond.zero} in this case.

\medskip

\noindent \textbf{Step 4. Case 2:} $\varphi'<0$ on $(0,R(a))$. Fix $\varrho_1\in(\varrho_0,R(a))$ and set
\begin{equation*}
	\delta:=1-m\varrho_1^{\sigma}V^{(m-1)/m}(\varrho_1)>0.
\end{equation*}
Then
\begin{equation*}
	1-mr^{\sigma}V^{(m-1)/m}(r)>\delta, \quad r\in(\varrho_1,R(a))
\end{equation*}
and integrating~\eqref{interm18} over $(\varrho_1,r)$ for $r\in(\varrho_1,R(a))$ gives
\begin{equation*}
\begin{split}
	r^{N-1}\varphi'(r)& = \varrho_1^{N-1} \varphi'(\varrho_1) + \int_{\varrho_1}^r \frac{s^{N-1}}{m} V^{(1-m)/m}(s) \varphi(s) \left[1-ms^{\sigma}V^{(m-1)/m}(s)\right]\,ds\\
	&\geq \varrho_1^{N-1}\varphi'(\varrho_1) + \frac{\delta\varrho_1^{N-1}}{m} \int_{\varrho_1}^r V^{(1-m)/m}(s) \varphi(s) \,ds.
\end{split}
\end{equation*}
The negativity of $\varphi'$ implies that $r^{N-1}\varphi'(r)\leq\varrho_1^{N-1}\varphi'(r)$ for any $r\in(\varrho_1,R(a))$ and we end up with
\begin{equation}\label{interm19}
	\varphi'(r)\geq\varphi'(\varrho_1)+\frac{\delta}{m}\int_{\varrho_1}^rV^{(1-m)/m}(s)\varphi(s)\,ds.
\end{equation}
On the one hand, the assumed positivity of $\varphi$ and negativity of $\varphi'$ on $(0,R(a))$, together with~\eqref{ds}, ensure that
\begin{equation*}
	\varphi''(r)=-\frac{N-1}{r}\varphi'(r)+\frac{1}{m} V^{(1-m)/m}(r)\varphi(r)\left[1-mr^{\sigma}V^{(m-1)/m}(r)\right]>0
\end{equation*}
for $r\in(\varrho_1,R(a))$. Therefore, there are $l_0\in[0,1)$ and \blm{$l_1\in [0,-\varphi'(\rho_1))$} such that
\begin{equation}\label{interm20}
\lim\limits_{r\to R(a)}\varphi(r)=l_0, \quad \lim\limits_{r\to R(a)}\varphi'(r)=-l_1, \quad \varphi(r)\sim l_0+l_1(R(a)-r), \ {\rm as} \ r\to R(a).
\end{equation}
On the other hand, according to the interface behavior~\eqref{expint} of $V$, \blm{there is $L(\rho_1)>C_I^{(m-1)/m}$ such that}
\begin{equation*}
	V^{(m-1)/m}(r)\leq \blm{L(\rho_1)} (R(a)-r)^2, \quad r\in(\varrho_1,R(a)),
\end{equation*}
hence
\begin{equation*}
V^{(1-m)/m}(r)\geq \blm{\frac{1}{L(\rho_1)(R(a)-r)^2},} \quad r\in(\varrho_1,R(a)).
\end{equation*}
We then deduce from~\eqref{interm19}, \blm{the assumed positivity of $\varphi$} and the negativity of $\varphi'$ \blm{on $(0,R(a))$} that
\begin{equation*}
	0\geq\varphi'(\varrho_1)+ \blm{\frac{\delta}{m L(\rho_1)}} \int_{\varrho_1}^r\frac{\varphi(s)}{(R(a)-s)^2}\,ds,
\end{equation*}
which in particular implies that
\begin{equation*}
	s\mapsto \frac{\varphi(s)}{(R(a)-s)^2}\in L^1(\varrho_1,R(a)),
\end{equation*}
in view of the positivity of $\varphi$. Consequently, we deduce from~\eqref{interm20} that $l_0=l_1=0$ and the claim \eqref{cond.zero} is established also in this case, completing the proof.
\end{proof}

The proof of the uniqueness of the zero of $\varphi$ is more involved and relies on the auxiliary functions we have introduced throughout this section.
\begin{proposition}\label{prop.uniqzero}
Let $a>0$ be such that $V(R(a))=V'(R(a))=0$. Then, there is a unique zero $r(a)\in(0,R(a))$ of $\varphi$ and
\begin{equation}\label{end.phi}
	\lim\limits_{r\to R(a)}\varphi(r) = -\infty, \quad \lim\limits_{r\to R(a)}\varphi'(r) = -\infty.
\end{equation}
\end{proposition}
\begin{proof}
We set
\begin{equation}\label{zero.phi}
	r(a):=\inf\{r\in(0,R(a)): \varphi(r;a)=0\} \in (0,R(a)),
\end{equation}
which is well defined according to Proposition~\ref{prop.existzero}. Then $\varphi(r(a))=0$ and, since $\varphi(0)=1$ by~\eqref{icphi}, necessarily $\varphi'(r(a))< 0$, the possibility that $\varphi'(r(a))=0$ being ruled out by the Cauchy-Lipschitz theorem and the fact that $\varphi$ is not identically equal to zero. We then deduce from~\eqref{deriv.xi} that $\xi'(r)>0$ for $r\in(0,r(a))$ and $\xi'(r(a))=0$. Since $\varphi'(r(a))<0$, we infer from~\eqref{deriv.xi} that
\begin{equation*}
	\xi''(r(a))=\frac{m-1}{m}r(a)^{N-1}V^{1/m}(r(a))\varphi'(r(a))<0,
\end{equation*}
hence $\xi'<0$ in a right neighborhood of $r(a)$. Moreover,~\eqref{exp.xizeta} ensures that $\xi(0)=0$, whence $\xi(r)>0$ for $r\in(0,r(a)]$. We split the rest of the proof into several steps.

\medskip

\noindent \textbf{Step 1. Positivity of the function $\xi$.} Our next goal is to prove that $\xi>0$ on $(0,R(a))$. To this end, we assume for contradiction that there exists $r_0\in(r(a),R(a))$ such that 
\begin{equation*}
	\xi(r)>0 \;\;\text{ for }\;\; r\in(0,r_0) \;\;\text{ and }\;\;  \xi(r_0)=0.
\end{equation*} 
Then $\xi'(r_0)\leq0$, which also implies that $\varphi(r_0)\leq0$ by~\eqref{deriv.xi}. If, moreover, $\xi'(r_0)=\varphi(r_0)=0$, then we deduce from the definition~\eqref{fxi} of $\xi$ that $V(r_0)\varphi'(r_0) =\xi(r_0)=0$, whence $\varphi(r_0)=\varphi'(r_0)=0$. Since $r_0\in(0,R(a))$, it follows by standard theory of differential equations that $\varphi\equiv0$ on $(0,R(a))$, which contradicts $\varphi(0)=1$. We have thus shown that
\begin{equation}\label{interm9}
	\varphi(r_0)<0, \quad \xi'(r_0)<0.
\end{equation}
It next follows from~\eqref{deriv.T}, along with the negativity of $\sigma$ and of $V'$ on $(0,R(a))$ and the definition of $r_0$, that
\begin{equation}\label{interm10}
	T'(r)<0, \quad r\in(0,r_0), \quad T'(r_0)=0.
\end{equation}
We further derive from~\eqref{fT}, \blm{\eqref{exp.T}} and~\eqref{interm10} that $T(r_0)=-\zeta(r_0)<0$. Since $\xi(r_0)=0$,~\eqref{fxi} entails that
\begin{equation*}
	V'(r_0)\varphi(r_0)=V(r_0)\varphi'(r_0)
\end{equation*}
and the definition~\eqref{fzeta} of $\zeta$ gives (taking into account~\eqref{fPQ})
\begin{equation*}
\begin{split}
	\zeta(r_0)&=\frac{\varphi(r_0)}{V(r_0)}\left[r_0^N\frac{V'(r_0)V(r_0)\varphi'(r_0)}{\varphi(r_0)}+r_0^{N}(r_0^{\sigma}V(r_0)-V^{1/m}(r_0))V(r_0)\right.\\
	&\hspace{3cm} \left. \phantom{\frac{\varphi(r_0)}{V(r_0)}} + (N-2)r_0^{N-1}V(r_0)V'(r_0)\right]\\
	&=\frac{\varphi(r_0)}{V(r_0)}\left[r_0^N(V')^2(r_0)+r_0^{N}(r_0^{\sigma}V(r_0)-V^{1/m}(r_0))V(r_0)\right.\\
	&\hspace{3cm} \left. \phantom{(V')^2(r_0)} + (N-2)r_0^{N-1}V(r_0)V'(r_0)\right]\\
	&=\frac{\varphi(r_0)}{V(r_0)}\left[P(r_0)+\frac{2m}{m+1}r_0^NV^{(m+1)/m}(r_0)-r_0^NV^{(m+1)/m}(r_0)\right],
\end{split}
\end{equation*}
hence
\begin{equation}\label{interm11}
	\zeta(r_0)=\frac{\varphi(r_0)}{V(r_0)}Q(r_0).
\end{equation}
Since $\varphi(r_0)<0$, $V(r_0)>0$ and $Q(r_0)>0$ by Lemma~\ref{lem.PQ}, we deduce from~\eqref{interm11} that $\zeta(r_0)<0$ and thus $T(r_0)=-\zeta(r_0)>0$, contradicting the negativity of $T(r_0)$ derived from~\eqref{interm10}. This contradiction proves that there is no point $r_0\in(r(a),R(a))$ such that $\xi(r_0)=0$ and, consequently, 
\begin{equation*}
	\xi>0 \;\;\text{ on }\;\; (0,R(a)).
\end{equation*} 
The positivity of $\xi$, along with~\eqref{deriv.T} and~\eqref{exp.T}, also guarantees that
\begin{equation}\label{interm13}
	T'(r)<0, \quad T(r)<0, \quad r\in(0,R(a)).
\end{equation}

\medskip

\noindent \textbf{Step 2. Non-existence of a second zero of $\varphi$.} Assume now for contradiction that there exists $\overline{r}(a)\in(r(a),R(a))$ such that $\varphi(r)<0$ for $r\in(r(a),\overline{r}(a))$ and $\varphi(\overline{r}(a))=0$. The latter also entails that $\varphi'(\overline{r}(a))>0$, since if
\begin{equation*}
	\varphi(\overline{r}(a))=\varphi'(\overline{r}(a))=0,
\end{equation*}
then $\varphi\equiv0$ on $(0,R(a))$, in contradiction with $\varphi(0)=1$. Since $\xi(r)>0$ for $r\in(0,R(a))$, we infer from~\eqref{fxi} that
\begin{equation}\label{interm12}
	V'(r)\varphi(r)>V(r)\varphi'(r), \quad r\in(0,R(a)).
\end{equation}
We evaluate~\eqref{interm12} at $r=\overline{r}(a)$ to obtain that $0>V(\overline{r}(a))\varphi'(\overline{r}(a))$, which contradicts the positivity of both functions $V$ and $\varphi'$ at $\overline{r}(a)$. We have just completed the proof of the uniqueness of a zero for $\varphi$ in $(0,R(a))$ and we also deduce that 
\begin{equation}
	\varphi(r)>0, \quad r\in (0,r(a)) \quad\text{ and }\quad\varphi(r)<0, \quad r\in(r(a),R(a)). \label{dphi}
\end{equation}

\medskip

\noindent \textbf{Step 3. Positivity of the function $\zeta$.} Setting
\begin{equation*}
	T_1(r):=V^{(m-1)/m}(r)-\frac{2}{\sigma+2}r^{-\sigma}, \quad r\in (0,R(a)),
\end{equation*}
it follows from the negativity of $\sigma$ and the definition of $R(a)$ that $T_1$ is decreasing on $(0,R(a))$. Moreover,
\begin{equation*}
	\lim\limits_{r\to0}T_1(r)=a^{(m-1)/m}>0, \quad \lim\limits_{r\to R(a)}T_1(r)=-\frac{2}{\sigma+2} R(a)^{-\sigma}<0.
\end{equation*}
Consequently, there is a unique $r_2(a)\in(0,R(a))$ such that $T_1(r)>0$ for $r\in(0,r_2(a))$ and $T_1(r)<0$ for $r\in(r_2(a),R(a))$. Since, by~\eqref{deriv.zeta},
\begin{equation*}
	\zeta'(r)=(\sigma+2)r^{\sigma+N-1}T_1(r)\varphi(r)V^{1/m}(r),
\end{equation*}
we find, taking into account the signs of $T_1$ and $\varphi$ on $(0,R(a))$, that
\begin{equation}\label{interm14}
\begin{split}
  &\zeta'(r)>0, \quad r\in \big(0,\min\{r(a),r_2(a)\}\big), \\
  &\zeta'(r)<0, \quad r\in \big(\min\{r(a),r_2(a)\},\max\{r(a),r_2(a)\}\big), \\
  &\zeta'(r)>0, \quad r\in \big(\max\{r(a),r_2(a)\},R(a)\big).
\end{split}
\end{equation}
Furthermore,~\eqref{fzeta}, \eqref{fT} and~\eqref{interm13}, along with the negativity of $V'(r(a))$ and $\varphi'(r(a))$, give
\begin{equation}\label{interm15}
\begin{split}
	\zeta(r(a)) & = r(a)^{N}V'(r(a))\varphi'(r(a))>0, \\
	\zeta(r_2(a)) & = \zeta(r_2(a)) - \frac{m(\sigma+2)}{m-1} r_2(a)^\sigma T_1(r_2(a)) \blm{\xi(r_2(a))} =-T(r_2(a))>0,
\end{split}
\end{equation}
Combining \blm{\eqref{exp.xizeta}, \eqref{interm14} and \eqref{interm15}} gives
\begin{equation}\label{pos.zeta}
	\zeta(r)\geq\zeta(\max\{r(a),r_2(a)\})>0,\quad  r\in\blm{\big(\min\{r(a),r_2(a)\},R(a)\big).}
\end{equation}

\medskip

\noindent \textbf{Step 4. Negativity of $\varphi$ and $\varphi'$ at $R(a)$.} We are only left to prove~\eqref{end.phi}. Let us first recall that, owing to~\eqref{ds}, \eqref{interm18} and~\eqref{dphi},
\begin{equation}\label{interm100}
\begin{split}
	& \big(r^{N-1}\varphi'(r)\big)' < 0 , \quad r\in \big(0,\min\{r(a),\rho_0\}\big), \\
	& \big(r^{N-1}\varphi'(r)\big)' > 0, \quad r\in \big(\min\{r(a),\rho_0\},\max\{r(a),\rho_0\}\big), \\
	& \big(r^{N-1}\varphi'(r)\big)' < 0, \quad r\in \big(\max\{r(a),\rho_0\},R(a)\big),
\end{split}
\end{equation}
and assume for contradiction that
\begin{equation}
	\inf\big\{\varphi'(r): r\in\big(\max\{r(a),\rho_0\},R(a)\big)\big\}\ge 0. \label{interm101}
\end{equation}
Then $\varphi$ is non-decreasing and negative on $\big(\max\{r(a),\rho_0\},R(a)\big)$ by~\eqref{dphi} and~\eqref{interm101}, so that there is $\varphi_I\le 0$ such that 
\begin{equation}
	\lim\limits_{r\to R(a)}\varphi(r) = \varphi_I \leq 0. \label{interm102}
\end{equation}	
If $\varphi_I<0$ we infer from~\eqref{expint} and~\eqref{interm18} that
\begin{equation*}
	 \big(r^{N-1}\varphi'(r)\big)' \sim \frac{R(a)^{N-1}}{m} \frac{C_I^{(1-m)/m}\varphi_I }{(R(a)-r)^{2}}, \quad {\rm as} \ r\to R(a), 
\end{equation*}	
hence, after integration,
\begin{equation*}
	r^{N-1}\varphi'(r) \sim \frac{R(a)^{N-1}}{m} \frac{C_I^{(1-m)/m}\varphi_I }{R(a)-r}, \quad {\rm as} \ r\to R(a).
\end{equation*}
Equivalently, since $R(a)\in (0,\infty)$ due to Lemma~\ref{lem.Ra},
\begin{equation*}
	0 \le \varphi'(r) \sim \frac{C_I^{(1-m)/m}\varphi_I }{m(R(a)-r)}, \quad {\rm as} \ r\to R(a),
\end{equation*}
which does not comply with the integrability of $\varphi$ as $r\to R(a)$ guaranteed by~\eqref{interm102}. Therefore, $\varphi_I=0$.

It next follows from~\eqref{fzeta}, the negativity of $V'$ and~\eqref{interm101} that
\begin{equation*}
	\zeta(r) \le r^N \big(r^\sigma V(r) - V^{1/m}(r)\big) \varphi(r) + (N-2) r^{N-1} (V'\varphi)(r)
\end{equation*}
for $r\in \big(\max\{r(a),\rho_0\},R(a)\big)$ and we may let $r\to R(a)$ in the above inequality and use~\eqref{pos.zeta} to conclude that
\begin{equation*}
	0 < \zeta(\max\{r(a),r_2(a)\}) \le \limsup_{r\to R(a)} \zeta(r) \le 0,
\end{equation*}
\blm{and a contradiction}. Therefore, \blm{the positivity property~\eqref{interm101} is not true, which implies, in view of~\eqref{interm100}, that} there are $(\delta_I,\varepsilon_I)\in (0,1)^2$ such that \blm{$R(a)-\delta_I>r(a)$ and}
\begin{equation}
	\varphi'(r) \le - \varepsilon_I, \quad r\in [R(a)-\delta_I,R(a))\blm{.} \label{interm103}
\end{equation}
\blm{An} immediate consequence of~\eqref{dphi} and~\eqref{interm103} is
\begin{equation*}
	\varphi(r) < \varphi(R(a)-\delta_I)<0, \quad r\in [R(a)-\delta_I,R(a)).
\end{equation*}
We then use the monotonicity~\eqref{interm103} of $\varphi$ on $(R(a)-\delta_I,R(a))$ and argue as above to show that
\begin{equation*}
	\lim_{r\to R(a)} \varphi(r) = \inf_{[R(a)-\delta_I,R(a))} \varphi = -\infty.
\end{equation*}
In particular, since $V(R(a))=0$ and $R(a)\in (0,\infty)$,
\begin{equation}
	\lim_{r\to R(a)} V^{(1-m)/m}(s) \varphi(s) \left[ 1 - m r^\sigma V^{(m-1)/m} (s)\right] = -\infty. \label{interm105}
\end{equation}
In addition, for $r\in (R(a)-\delta_I,R(a))$, it follows \blm{from~\eqref{interm18}} and~\eqref{interm103} that
\begin{align*}
	r^{N-1} \varphi'(r) & = (R(a)-\delta_I)^{N-1} \varphi'(R(a)-\delta_I) \\
	& \quad + \frac{1}{m} \int_{R(a)-\delta_I}^r s^{N-1} V^{(1-m)/m}(s) \varphi(s) \left[ 1 - m r^\sigma V^{(m-1)/m} (s)\right]\,ds \\
	& \le \frac{1}{m} \int_{R(a)-\delta_I}^r s^{N-1} V^{(1-m)/m}(s) \varphi(s) \left[ 1 - m r^\sigma V^{(m-1)/m} (s)\right]\,ds
\end{align*}
and we let $r\to R(a)$ in the above inequality to deduce that $r^{N-1} \varphi'(r)\longrightarrow -\infty$ as $r\to R(a)$. Thus, both statements in~\eqref{end.phi} hold true and the proof is complete.
\end{proof}

\section{Proof of Theorem \ref{th.uniq}}\label{sec.proof}

We are now in a position to complete the classification with respect to $a$ of the behavior of $V(\cdot;a)$ \blm{at $R(a)$} and thus the proof of Theorem~\ref{th.uniq}. The energy $E(\cdot;a)$ introduced in Lemma~\ref{lem.energy} plays an important role in the next lines, and we recall it here for the reader's convenience
\begin{equation}\label{energybis}
	E(r;a)=\frac{1}{2}(V')^2(r;a) \blm{+\frac{1}{2}r^{\sigma}V^2(r;a) - \frac{m}{m+1}V^{(m+1)/m}(r;a)}, \quad r\in (0,R(a)].
\end{equation}
From the definition of $R(a)$ and Lemma~\ref{lem.Ra}, we know that $V(\cdot;a)$ is positive and decreasing on $(0,R(a))$ with $R(a)\in (0,\infty)$ and
\begin{equation}\label{limits}
	V(R(a);a)\in [0,a), \quad V'(R(a);a)\in(-\infty,0].
\end{equation}
We can thus introduce the following three sets
\begin{equation}\label{sets}
\begin{split}
	&\mathcal{A}:=\left\{a\in(0,\infty)\ \blm{:} \ V(R(a);a)=0 \ {\rm and} \ V'(R(a);a)<0 \right\},\\
	&\mathcal{B}:=\left\{a\in(0,\infty)\ : V(R(a);a)=V'(R(a);a)=0\right\},\\
	&\mathcal{C}:=\left\{a\in(0,\infty)\ : V(R(a);a)>0 \right\}.
\end{split}
\end{equation}
We infer from~\eqref{limits} and the definitions~\eqref{sets} that the three sets $\mathcal{A}$, $\mathcal{B}$ and $\mathcal{C}$ form a partition of $(0,\infty)$. They can be also fully characterized in terms of the energy~\eqref{energybis}, as follows.
\begin{proposition}\label{prop.energy}
The following equivalent classification holds true:
\begin{equation*}
\begin{split}
&\mathcal{A}:=\left\{a\in(0,\infty)\ :\ E(R(a);a)>0 \right\},\\
&\mathcal{B}:=\left\{a\in(0,\infty)\ :\ E(R(a);a)=0 \right\},\\
&\mathcal{C}:=\left\{a\in(0,\infty)\ : \ E(R(a);a)<0 \right\}.
\end{split}
\end{equation*}
Moreover, if $a\in\mathcal{C}$, then $V''(R(a);a)>0$ and $V'(\cdot;a)$ is positive in a right neighborhood of $R(a)$.
\end{proposition}
\begin{proof}
It is obvious from~\eqref{energybis} and the definitions of $\mathcal{A}$ and $\mathcal{B}$ that 
\begin{equation*}
	E(R(a);a)>0 \;\;\text{ for }\;\; a\in\mathcal{A} \;\;\;\text{ and }\;\;\; E(R(a);a)=0 \;\;\text{ for }\;\; a\in\mathcal{B}. 
\end{equation*}
Finally, pick $a\in\mathcal{C}$. Owing to the definition and finiteness of $R(a)$, there holds
\begin{equation}
	V'(R(a);a) = 0. \label{interm104}
\end{equation}
\blm{As in the proof of Lemma~\ref{lem.Ra}, consider} the function
\begin{equation*}
	\blm{T_0}(r;a):=1-r^{\sigma}V^{(m-1)/m}(r;a), \quad r\in (0,R(a)).
\end{equation*}
We observe that $\blm{T_0}(r;a)\to-\infty$ as $r\to 0$, and $\blm{T_0(\cdot;a)}$ is increasing on $(0,R(a))$ in view of the negativity of $\sigma$ and $V'(\cdot;a)$ and the fact that $m>1$. Moreover, we deduce \blm{from~\eqref{CPV}} that 
\begin{equation*}
	r^{N-1}V'(r;a)=\int_0^r s^{N-1}V^{1/m}(s;a)\blm{T_0}(s;a)\,ds, \quad r\in (0,R(a)).
\end{equation*}
Letting $r\to R(a)$ in the above identity and using~\eqref{interm104}, we find 
\begin{equation}\label{interm23}
	\int_0^{R(a)}s^{N-1}V^{1/m}(s;a)\blm{T_0}(s;a)\,ds=0.
\end{equation}
The positivity of $V(\cdot;a)$, along with~\eqref{interm23}, then implies that there is $s_0(a)\in(0,R(a))$ such that $\blm{T_0}(s_0(a);a)=0$. \blm{Therefore}, $\blm{T_0}(R(a);a)>0$, which ensures that
\begin{equation}\label{interm24}
	R(a)^{\sigma}V^{(m-1)/m}(R(a);a)<1,
\end{equation}
and
\begin{equation}
\begin{split}
	V''(R(a);a) & = - \frac{N-1}{R(a)} \blm{V'}(R(a);a) + V^{1/m}(R(a);a) \blm{T_0}(R(a);a) \\
	& = V^{1/m}(R(a);a) \blm{T_0}(R(a);a)  > 0. 
\end{split} \label{pam}
\end{equation}
We finally infer from~\eqref{energybis}, \eqref{interm104} and~\eqref{interm24} that 
\begin{equation*}
\begin{split}
	E(R(a);a) & = \frac{1}{2} V(R(a);a)^{(m+1)/m} \left[-\frac{2m}{m+1}+R(a)^{\sigma}V^{(m-1)/m}(R(a);a)\right]\\
	&<-\frac{m-1}{2(m+1)}V(R(a);a)^{(m+1)/m}<0,
\end{split}
\end{equation*}
concluding the proof.
\end{proof}

We next show that the sets $\mathcal{A}$ and $\mathcal{C}$ are open, a decisive role being played by Proposition~\ref{prop.energy} and the extension of $V(\cdot;a)$ beyond $R(a)$ in this step. 
\begin{proposition}\label{prop.open}
The sets $\mathcal{A}$ and $\mathcal{C}$ are open.
\end{proposition}
\begin{proof}
We first prove that $\mathcal{C}$ is open. Let $a_0\in\mathcal{C}$. It follows from the definition of $\mathcal{C}$ and Proposition~\ref{prop.energy} that there are $\delta\in (0,R(a_0)/2)$ and $\varepsilon\in (0,1)$ such that
\begin{equation}
\begin{split}
	\inf_{r\in [0,R(a_0)-\delta]} V(r;a_0) \ge 2\varepsilon, \quad & \quad \sup_{r\in [\delta,R(a_0)-\delta]} V'(r;a_0) \le - 2\varepsilon, \\
	2^{1/\sigma} Y(a_0) \ge 2\delta, \quad & \quad  E(R(a_0)-2\delta;a_0) \le - 2 \varepsilon,
\end{split}\label{cd01}
\end{equation}
where the function $Y(a)$ is the unique solution to the equation
\begin{equation*}
	Y(a)^{-m\sigma/(m-1)}+\frac{aY(a)^{\sigma+2}}{(\sigma+2)(\sigma+N)}=a
\end{equation*}
defined in Lemma~\ref{lem.propV}. By the continuous dependence established in Proposition~\ref{prop.stab} and the continuity of the function $Y$, there are $0<a_0/2\le a_{\min}<a_0<a_{\max}\le 2a_0$ such that
\begin{equation*}
	\inf_{r\in [0,R(a_0)-\delta]} V(r;a) \ge \varepsilon, \quad \sup_{r\in [\delta,R(a_0)-\delta]} V'(r;a) \le - \varepsilon, \quad 2^{1/\sigma} Y(a) \ge \delta, \quad a\in (a_{\min},a_{\max}). 
\end{equation*}
Combining the above properties with Corollary~\ref{cor.propV} ensures that $R(a)>R(a_0)-\delta$ and we use again the continuous dependence established in Proposition~\ref{prop.stab} to conclude that $E(R(a_0)-2\delta;a)\le -\varepsilon$ for $a\in (a_{\min},a_{\max})$, after possibly reducing the range of this interval. It then follows from the monotonicity of $E(\cdot;a)$ that 
\begin{equation*}
	E(R(a);a) \le E(R(a_0)-2\delta;a)\le -\varepsilon < 0, \quad a\in (a_{\min},a_{\max}),
\end{equation*}
and Proposition~\ref{prop.energy} implies that $(a_{\min},a_{\max})\subset \mathcal{C}$.

We next prove that $\mathcal{A}$ is open. Pick $a_0\in\mathcal{A}$. According to the definition of $\mathcal{A}$, there are $\delta\in (0,R(a_0)/2)$ and $\varepsilon\in (0,1)$ such that
\begin{equation*}
\begin{split}
	\inf\limits_{r\in[0,R(a_0)-\delta]}V(r;a_0)\geq 2\varepsilon, \quad & \quad V(R(a_0)+\delta;a_0) \le - 2\varepsilon, \\
	\sup\limits_{r\in[\delta,R(a_0)+\delta]}V'(r;a_0) \le -2\varepsilon, \quad & \quad 2^{1/\sigma} Y(a_0) \ge 2\delta, \quad R(a_0) + 2\delta \le S(a_0),
\end{split}
\end{equation*}
 the last property being a consequence of Proposition~\ref{prop.wpv} which guarantees that $R(a_0)<S(a_0)$. As before, we infer from Corollary~\ref{cor.propV}, Proposition~\ref{prop.stab} and the continuity of the function $\blm{Y}$ defined in Lemma~\ref{lem.propV} that there are $0<a_0/2\le a_{\min}<a_0<a_{\max}\le 2a_0$ such that
\begin{equation*}
\begin{split}
	\inf\limits_{r\in[0,R(a_0)-\delta]}V(r;a)\geq \varepsilon, \quad & \quad V(R(a_0)+\delta;a) \le - \varepsilon, \\
	\sup\limits_{r\in[\delta,R(a_0)+\delta]}V'(r;a) \le -\varepsilon, \quad & \quad V'(r;a) < 0, \quad r\in (0,\delta).
\end{split}
\end{equation*}
Therefore, $R(a)\in (R(a_0)-\delta,R(a_0)+\delta)$ with $V'(R(a);a)\le -\varepsilon$ for any $a\in(a_{\min},a_{\max})$ and we have shown that $(a_{\min},a_{\max})\subset \mathcal{A}$; that is, $\mathcal{A}$ is open and the proof is complete.
\end{proof}
After all these preparations, we are now in a position to conclude the proof of Theorem~\ref{th.uniq}, which follows similar ideas as in the proof of the analogous result in \cite{CEF96}.

\begin{proof}[Proof of Theorem~\ref{th.uniq}]
We first recall that, by \cite[Theorems~1.1 and~1.2]{IL25}, the set $\mathcal{B}$ is non-empty. We thus consider $a_0\in\mathcal{B}$ and let $\varepsilon>0$ be such that 
\begin{equation}\label{interm25}
	(2\varepsilon)^{(1-m)/m}>m\left(\frac{R(a_0)}{2}\right)^{\sigma}.
\end{equation}
We infer from Propositions~\ref{prop.existzero} and~\ref{prop.uniqzero} that there are $r_0\in(R(a_0)/2,R(a_0))$ and $\delta_0\in (0,r_0)$ such that 
\begin{equation}
\begin{split}
	\inf_{r\in[0,r_0]} V(r;a_0)\ge 2\varepsilon, \quad & \quad 0 \leq V(r;a_0) \leq 2\varepsilon, \quad r\in[r_0,R(a_0)], \\
	\sup_{r\in[\delta_0,r_0]} V'(r;a_0) \le -2\varepsilon, \quad & \quad 2^{1/\sigma} Y(a_0)\ge 2 \delta_0,
\end{split} \label{z00}
\end{equation}
and 
\begin{equation*}
	\varphi(r_0;a_0) \le -2\varepsilon, \quad \varphi'(r_0;a_0) \le -2\varepsilon.
\end{equation*}
By Corollary~\ref{cor.propV} and Propositions~\ref{prop.stab} and~\ref{prop.phi}, there is $\delta>0$ such that, for any $a\in(a_0-\delta,a_0+\delta)$,
\begin{equation}
\begin{split}
	\inf_{r\in[0,r_0]}V(r;a) \ge \varepsilon, \quad & \quad -\varepsilon \leq V(r;a)\leq \varepsilon, \quad r\in[r_0,R(a_0)], \\
	\sup_{r\in[\delta_0,r_0]} V'(r;a) \le -\varepsilon, \quad & \quad V'(r;a)<0, \quad r\in (0,\delta_0),
\end{split}\label{z0}
\end{equation}
and
\begin{equation}
	\varphi(r_0;a) \le -\varepsilon, \quad \varphi'(r_0;a) \le -\varepsilon. \label{z1}
\end{equation}
We also note for further use that~\eqref{interm25} ensures that 
\begin{equation}
	z\mapsto r^{\sigma} z - z^{1/m} \;\text{ is decreasing on }\; (0,2\varepsilon), \quad \text{ for all }\; r\in[r_0,R(a_0)]. \label{z2}
\end{equation}
Indeed, for $r\in [r_0,R(a_0)]\subset [R(a_0)/2,R(a_0)]$ and $z\in (0,2\varepsilon)$,
\begin{equation*}
	\frac{d}{dz}\big(r^{\sigma}z-z^{1/m}\big) = r^{\sigma} - \frac{1}{m} z^{(1-m)/m} < \left( \frac{R(a_0)}{2} \right)^{\sigma} - \frac{1}{m} (2\varepsilon)^{(1-m)/m} <0.
\end{equation*}
We split the rest of the analysis into two steps.

\medskip 

\noindent \textbf{Step~1.} Let $a\in(a_0,a_0+\delta)$. Then we deduce from~\eqref{z1} that
\begin{equation*}
	V(r_0;a)<V(r_0;a_0), \quad V'(r_0;a)<V'(r_0;a_0).
\end{equation*}
In addition, $V(r_0;a)\ge \varepsilon>0$ by~\eqref{z0} and $V(R(a_0);a_0)=0$, so that, in view of~\eqref{z00} and~\eqref{z2}, we are in a position to apply~Proposition~\ref{prop.comp}~(b) with $(r_1,r_2) = (r_0,R(a_0))$ and $(h,g)=\big(V(\cdot;a),V(\cdot;a_0)\big)$ to conclude that there is $\rho_a\in (r_0,R(a_0))$ such that $V(\cdot;a)>0$ on $[r_0,\rho_a)$, $V(\rho_a;a)=0$ and $V'(\rho_a;a)<0$; that is, $R(a)=\rho_a$ and $a\in \mathcal{A}$, \blm{owing to~\eqref{z0}}. We have just proved that $(a_0,a_0+\delta)\subseteq\mathcal{A}$.

\medskip 

\noindent \textbf{Step~2.} Let $a\in(a_0-\delta,a_0)$. Then
\begin{equation}
	V(r_0;a_0)<V(r_0;a), \quad V'(r_0;a_0)<V'(r_0;a) \label{z3}
\end{equation}
by~\eqref{z1}. Assume for contradiction that $R(a)\ge R(a_0)$. Then $V(\cdot;a)$ is non-negative on $[r_0,R(a_0)]$ and we infer from~\eqref{z0}, \eqref{z2}, \eqref{z3} and Proposition~\ref{prop.comp}~(a) with $(r_1,r_2)=(r_0,R(a_0))$ and $(h,g)=\big(V(\cdot;a_0),V(\cdot;a)\big)$ that 
\begin{equation*}
	0 \le V(r;a_0) < V(r;a), \quad r\in [r_0,R(a_0)], \quad\quad 0 = V'(R(a_0);a_0) < V'(R(a_0);a)\blm{.}
\end{equation*}
Since $V'(r_0;a)\le -\varepsilon<0$ by~\eqref{z0}, Bolzano's theorem ensures the existence of $\rho\in (r_0,R(a_0))$ such that $V'(\rho;a)=0$. Consequently, $R(a)\le \rho<R(a_0)$ and a contradiction. 

We have thus shown that $R(a)<R(a_0)$. Since $V(\cdot;a)$ is non-negative on $[r_0,R(a)]$, we apply once more Proposition~\ref{prop.comp}~(a), still with $(h,g)=\big(V(\cdot;a_0),V(\cdot;a)\big)$ but now on $(r_1,r_2)=(r_0,R(a))$, to obtain
\begin{equation*}
	0 \le V(r;a_0) < V(r;a), \quad r\in [r_0,R(a)], \quad\quad V'(R(a);a_0) < V'(R(a);a).
\end{equation*}
In particular, $V(R(a);a)>0$, so that $a\in\mathcal{C}$ and we have just proved that $(a_0-\delta,a_0)\subseteq\mathcal{C}$.

\medskip 

\noindent \textbf{Step~3. End of the proof.} Define $\bar{a} :=\sup\{a>a_0\ :\ (a_0,a)\subset \mathcal{A}\}$. According to \textbf{Step~1}, $\bar{a}>a_0$. Assume for contradiction that $\bar{a}<\infty$. Then $(a_0, \bar{a})\subset\mathcal{A}$ by the definition of $\bar{a}$ and $\bar{a}\not\in\mathcal{A}$ since $\mathcal{A}$ is open by Proposition~\ref{prop.open}. Similarly, $\bar{a}\not\in\mathcal{C}$, again using the fact that $\mathcal{C}$ is open by Proposition~\ref{prop.open}. Therefore, $\bar{a}\in\mathcal{B}$ but the inclusion $(a_0,\bar{a})\subset\mathcal{A}$ contradicts \textbf{Step~2}. Consequently, $\bar{a}=\infty$ and $(a_0,\infty)\subset \mathcal{A}$. A similar argument shows that $(0,a_0)\subset \mathcal{C}$ and entails that $\mathcal{B}=\{a_0\}$ with $\mathcal{A}=(a_0,\infty)$ and $\mathcal{C}=(0,a_0)$, completing the proof. 
\end{proof}

\section{Proof of Theorem \ref{th.ckn01}}\label{sec.CKN}

We complete the paper with the proof of Theorem~\ref{th.ckn01}, which is a straightforward consequence of Theorem~\ref{th.uniq} and the following result.
\begin{proposition}\label{prop.ckn01}
Let $w\in \mathcal{X}$, $w\ne 0$, be a non-negative minimizer of~\eqref{ckni}; that is,
\begin{equation}
	w\ge 0 \;\text{ and }\; \mathcal{S}(w) = \mathcal{S}_*. \label{ckn03}
\end{equation} 
Then $w$ is a radially symmetric function with non-increasing profile and there are $(\zeta,\eta)\in (0,\infty)^2$ such that $x\mapsto \zeta w(\eta x)$ solves the Cauchy problem \eqref{CPV}.
\end{proposition}
Similarly as in \cite{IL25,ZZ15}, an intermediate step in the proof is to show that minimizers of~\eqref{ckni} are closely connected to minimizers of a constrained variational problem. Specifically, adopting for simplicity the notation
\begin{equation*}
	p:=\frac{m+1}{m}\in(1,2), \quad \lambda:=\frac{2m}{m+1}=\frac{2}{p},
\end{equation*}
and introducing
\begin{equation*}
	\mathscr{A}:= \left\{ v\in\mathcal{X}\ :\ N_\sigma(v)=1 \right\}
\end{equation*}
and the functional
\begin{equation}
	\mathcal{I}(v) := \|\nabla v\|_2^2 + \lambda \|v\|_{\blm{p}}^{\blm{p}}, \quad v\in\mathcal{X}, \label{ckn04}
\end{equation}
we first prove the following auxiliary result, see \cite[Theorem~3.1]{ZZ15}.

\begin{lemma}\label{lem.ckn02}
For $v\in\mathscr{A}$ and $\xi>0$, we define $v_\xi(x) := \xi^{(\sigma+N)/2} v(\xi x)$, $x\in\mathbb{R}^N$. Then,
\begin{equation*}
	\inf_{\xi>0} \mathcal{I}(v_\xi) = \mathcal{I}(v_{\xi_v}) = \kappa \big[\mathcal{S}(v)\big]^\nu
\end{equation*}
with
\begin{align*}
	\alpha & := \sigma+2>0, \quad \beta := \frac{N(m-1)-\sigma(m+1)}{2m}>0, \\ 
	\nu & := \frac{N(m-1)+2(m+1)}{(\sigma+N)(m-1)+4m}, \\ 
	\kappa & := \lambda^{\alpha/(\alpha+\beta)} (\alpha+\beta) \alpha^{-\alpha/(\alpha+\beta)} \beta^{-\beta/(\alpha+\beta)}
\end{align*}
and
\begin{equation*}
	\xi_v := \left( \frac{\beta \lambda \|v\|_p^p}{\alpha \|\nabla v\|_2^2} \right)^{1/(\alpha+\beta)}.
\end{equation*}
\end{lemma}

\begin{proof}
We sketch the proof for the sake of completeness. We first note that, for $v\in\mathscr{A}$ and $\xi>0$,
\begin{equation*}
	N_\sigma(v_\xi) = N_\sigma(v) = 1, \quad	\|\nabla v_\xi\|_2^2 = \xi^\alpha \|\nabla v\|_2^2, \quad \|v_\xi\|_p^p = \xi^{-\beta} \|v\|_p^p, 
\end{equation*} 
so that $v_\xi\in\mathscr{A}$ with
\begin{equation*}
	\mathcal{I}(v_\xi) = \xi^\alpha \|\nabla v\|_2^2 + \lambda \xi^{-\beta} \|v\|_p^p.
\end{equation*}
Therefore, \blm{recalling the definition~\eqref{ckne} of $a$,}
\begin{align*}
	\inf_{\xi>0} \mathcal{I}(v_\xi) & = \mathcal{I}(v_{\xi_v}) = \kappa \|\nabla v\|_2^{2\beta/(\alpha+\beta)} \|v\|_p^{p\alpha/(\alpha+\beta)} \\
	& = \kappa \left[ \|\nabla v\|_2^{2a} \|v\|_p^{2(1-a)} \right]^\nu,
\end{align*}
and the proof is complete.
\end{proof}

\begin{lemma}\label{lem.ckn03}
In the same conditions and notation as above, we have
\begin{equation}
	\mathcal{I}_* := \inf_{v\in\mathscr{A}} \mathcal{I}(v) = \kappa \mathcal{S}_*^\nu. \label{ckn05}
\end{equation}
\end{lemma}

\begin{proof}
On the one hand, we infer from Lemma~\ref{lem.ckn02} that, for $v\in\mathscr{A}$, 
\begin{equation*}
	\mathcal{I}(v) \ge \inf_{\xi>0} \mathcal{I}(v_\xi) = \kappa \big[\mathcal{S}(v)\big]^\nu \ge \kappa \mathcal{S}_*^\nu.
\end{equation*}
Therefore,
\begin{equation*}
	\mathcal{I}_* \ge \kappa \mathcal{S}_*^\nu.
\end{equation*}
On the other hand, since $\{v_\xi\ :\ \xi>0\}\subset \mathscr{A}$ for $v\in\mathscr{A}$, 
\begin{equation*}
	\mathcal{I}_* \le \inf_{\xi>0} \mathcal{I}(v_\xi) = \kappa \big[\mathcal{S}(v)\big]^\nu, 
\end{equation*}
from which we deduce that
\begin{equation*}
	\mathcal{I}_* \le \kappa \left[ \inf_{v\in\mathcal{A}} \mathcal{S}(v) \right]^\nu.
\end{equation*}
Now, given $v\in\mathcal{X}$, $v\ne 0$, the function $\bar{v}=v/N_\sigma(v)$ belongs to $\mathscr{A}$ with $\mathcal{S}(\bar{v}) = \mathcal{S}(v)$, hence
\begin{equation*}
	\inf_{v\in\mathscr{A}} \mathcal{S}(v) = \mathcal{S}_*. 
\end{equation*}
Gathering the outcome of the previous computations readily leads to~\eqref{ckn05}.
\end{proof}

\begin{lemma}\label{lem.ckn04}
Let $w\in \mathscr{A}$ be a non-negative minimizer of $\mathcal{I}$ in $\mathscr{A}$; that is,
\begin{equation}
	w\ge 0 \;\text{ and }\; \mathcal{I}(w) = \mathcal{I}_*. \label{ckn06}
\end{equation} 
Then $w$ is radially symmetric with non-increasing profile.
\end{lemma}

\begin{proof}
We further develop the proof of \cite[Lemma~2.2]{IL25}. Let $w^*$ be the symmetric-decreasing rearrangement of $w$. Owing to the pointwise identity $\big(w^2)^*=\big(w^*\big)^2$ \cite[Eq.~(3.8)]{Li77}, the monotonicity of $x\mapsto |x|^\sigma$ and the Hardy-Littlewood inequality \cite[3.4~Theorem]{LL01}, we have
\begin{equation}
	1 = N_\sigma^2(w) \le \blm{\int_{\mathbb{R}^N} |x|^\sigma \big(w^2\big)^*(x)\,dx} = N_\sigma^2\big(w^*\big). \label{ckn07}
\end{equation}
Since $W:=w^*/N_\sigma\big(w^*\big)$ belongs to $\mathscr{A}$, we infer from the equimeasurability property of the symmetric-decreasing rearrangement, the fact that $p<2$, Polya's inequality \cite[7.17~Lemma]{LL01}, \eqref{ckn06} and~\eqref{ckn07} that
\begin{align*}
	\mathcal{I}_* & = \mathcal{I}(w) \le \mathcal{I}(W) = \frac{\|\nabla w^*\|_2^2}{N_\sigma^2\big(w^*\big)} + \lambda \frac{\big\|w^*\big\|_p^p}{N_\sigma^p\big(w^*\big)} \\
	& \le \frac{\|\nabla w\|_2^2 + \lambda \|w\|_p^p}{N_\sigma^p\big(w^*\big)} = \frac{\mathcal{I}(w)}{N_\sigma^p\big(w^*\big)} = \frac{\mathcal{I}_*}{N_\sigma^p\big(w^*\big)}.
\end{align*}
Recalling that $\mathcal{I}_*>0$ according to the definition of $\mathscr{A}$, we conclude that $N_\sigma^p\big(w^*\big)\le 1$, which we combine with~\eqref{ckn07} to obtain
\begin{equation}
	N_\sigma\big(w^*\big) = N_\sigma(w) = 1. \label{ckn08}
\end{equation}
Now, since $x\mapsto |x|^\sigma$ is strictly symmetric-decreasing in the sense of \cite[3.3~(ii)]{LL01}, it follows from~\eqref{ckn08} and \cite[3.4~Theorem]{LL01} that $w=w^*$, which completes the proof.
\end{proof}

We can now complete the proof of both Proposition~\ref{prop.ckn01} \blm{and} Theorem~\ref{th.ckn01}.

\begin{proof}[Proof of Proposition~\ref{prop.ckn01} and Theorem~\ref{th.ckn01}]
Let $w\in \mathcal{X}$, $w\ne 0$, be a non-negative minimizer of~\eqref{ckni}. Then $W=w/N_\sigma(w)\in\mathscr{A}$ also satisfies~\eqref{ckn03}. By Lemma~\ref{lem.ckn02} and Lemma~\ref{lem.ckn03}, 
\begin{equation*}
	\mathcal{I}(W_{\xi_W}) = \kappa \big[\mathcal{S}(W)\big]^\nu = \kappa \mathcal{S}_*^\nu = \mathcal{I}_*,
\end{equation*}
recalling that $\xi_W$ is defined in Lemma~\ref{lem.ckn02}. Consequently, $W_{\xi_W}$ is a minimizer of $\mathcal{I}$ in $\mathscr{A}$ and Lemma~\ref{lem.ckn04} implies that $W_{\xi_W}$ is radially symmetric with non-increasing profile. Since
\begin{equation*}
	w(x) = \frac{N_\sigma(w)}{\xi_W^{(\sigma+N)/2}} W_{\xi_W}\left( \frac{x}{\xi_W} \right), \quad x\in\mathbb{R}^N,
\end{equation*}
we readily conclude that $w$ is also radially symmetric with non-increasing profile.
	
Furthermore, arguing as in the proof of \cite[Proposition~2.3]{IL25}, we deduce from the minimizing property of $W_{\xi_W}$ that it is a variational solution to
\begin{equation*}
	- \Delta W_{\xi_W}(x) + W_{\xi_W}^{1/m}(x) = \rho |x|^\sigma W_{\xi_W}(x), \quad x\in\mathbb{R}^N,
\end{equation*}
with $\rho := \|\nabla W_{\xi_W}\|_2^2 + \|W_{\xi_W}\|_p^p$. Setting 
\begin{equation*}
	\bar{w}(x) = \rho^{2m/[(m-1)(\sigma+2)]} W_{\xi_W}\big(\rho^{-1/(\sigma+2)}x\big), \quad x\in\mathbb{R}^N, 
\end{equation*}
we conclude that $\bar{w}$ is a variational solution to 
\begin{equation*}
	- \Delta \bar{w}(x) + \bar{w}^{1/m}(x) = |x|^\sigma \bar{w}(x), \quad x\in\mathbb{R}^N,
\end{equation*}
and
\begin{equation*}
	w(x) = \frac{N_\sigma(w)}{\rho^{2m/[(m-1)(\sigma+2)]} \xi_W^{(\sigma+N)/2}} \bar{w}\left( \frac{x \rho^{1/(\sigma+2)}}{\xi_W} \right), \quad x\in\mathbb{R}^N,
\end{equation*}
which proves the last statement of Theorem~\ref{th.ckn01} after setting 
\begin{equation*}
	\zeta:= \frac{N_\sigma(w)}{\rho^{2m/[(m-1)(\sigma+2)]} \xi_W^{(\sigma+N)/2}}, \quad \eta := \frac{\rho^{1/(\sigma+2)}}{\xi_W},
\end{equation*}
and noticing that $\bar{w}$ is also radially symmetric with non-increasing profile. The proof of Theorem~\ref{th.ckn01} is completed by the uniqueness established in Theorem~\ref{th.uniq}, which allows us to identify $\bar{w}=V(\cdot;a_*)$.
\end{proof}

\appendix
\section{Well-posedness for a linear initial value problem with singular coefficients}\label{app.A}
This first technical section is dedicated to establishing the well-posedness of a second-order linear ordinary differential equation with coefficients featuring a singularity at the origin, which will be then used repeatedly to obtain well-posedness and differentiability properties for solutions to Eqs.~\eqref{eq1.rad} and~\eqref{eq1.phi}. 

For $R>0$, let $U\in L^{1}_{\rm loc}([0,R))$ and $F\in L^{1}_{\rm loc}([0,R))$ be two locally integrable functions and consider the following linear initial value problem
\begin{subequations}\label{singoper}
\begin{equation}\label{L0}
	h''(r)+\frac{N-1}{r}h'(r)+r^{\sigma}h(r)-U(r)h(r)=F(r), \quad r\in(0,R),
\end{equation}
\begin{equation}\label{ic}
	h(0)=h_0\in\mathbb{R}, \quad \lim\limits_{r\to0}r^{N-1}h'(r)=0.
\end{equation}
\end{subequations}
Observe that~\eqref{L0} features a possibly non-integrable but algebraic singularity, as well as an integrable singularity, both as $r\to 0$. 

We first prove uniqueness and continuous dependence for the Cauchy problem~\eqref{singoper} for solutions in $C^2((0,R))\cap C([0,R))$, the building block of the proof being the estimate which we derive now.
\begin{lemma}\label{lem.uniq}
Let $h\in C^2((0,R))\cap C([0,R))$ be a solution to the Cauchy problem~\eqref{singoper}. Then, for any $R_0\in(0,R)$, there is a positive constant \blm{$C_0(R_0)$, depending only on $R_0$, such} that, for all $r\in (0,R_0)$,
\begin{eqnarray}
	|h(r)| & \leq & \blm{C_0(R_0) \big(|h_0| + \|F\|_{L^1(0,R_0)}\big)} \exp\left\{ C_0(R_0) \left(\frac{r^{\sigma_N}}{\sigma_N} + r^{\sigma_N-\sigma-1}\|U\|_{L^1(0,r)}\right) \right\}, \label{esth} \\
	r^{N-1} |h'(r)| & \leq & r^{N-1} \|F\|_{L^1(0,r)} + \left( \frac{r^{\sigma+N}}{\sigma+N} + r^{N-1}\|U\|_{L^1(0,r)}\right) \|h\|_{L^\infty(0,r)}, \label{esthp} 
\end{eqnarray}
where
\begin{equation*}
\sigma_N:=
\begin{cases}
    \sigma+1, & \mbox{if } N=1, \\
    \sigma+2-\epsilon, & \mbox{if } N=2, \\
    \sigma+2, & \mbox{if } N\geq3,
\end{cases}
\end{equation*}
for a fixed but arbitrary $\epsilon\in (0,\sigma+2)\cap (0,1)$ in dimension $N=2$. If, moreover,
\begin{equation}\label{global}
	U\in L^1(0,R), \quad F\in L^{1}(0,R),
\end{equation}
the estimates~\eqref{esth} and~\eqref{esthp} hold true for $R_0=R$ and any $r\in[0,R]$.
\end{lemma}

Note that, in all the three cases, $\sigma_N>0$ and also $\sigma_N\geq\sigma+1$.
\begin{proof}
We infer from~\eqref{L0} that
\begin{equation*}
\begin{split}
	\frac{d}{dr}\left(r^{N-1}h'(r)\right)&=r^{N-1}\left(h''(r)+\frac{N-1}{r}h'(r)\right)\\
	&=-r^{N+\sigma-1}h(r)+r^{N-1}U(r)h(r)+r^{N-1}F(r),
\end{split}
\end{equation*}
and we integrate over $(0,r)$ and take into account the second condition in~\eqref{ic} to obtain that
\begin{equation}\label{interm1}
	h'(r)=\int_0^r\left(\frac{s}{r}\right)^{N-1}\left(-s^{\sigma}h(s)+U(s)h(s)+F(s)\right)\,ds.
\end{equation}
On the one hand, the estimate~\eqref{esthp} readily follows from~\eqref{interm1} in view of the positivity of  $\sigma+N$. On the other hand, we integrate again the identity~\eqref{interm1} over $(0,r)$ and, since $h(0)=h_0$, we find\blm{,} by Fubini's theorem\blm{,}
\begin{equation}\label{interm2}
\begin{split}
h(r)&=h_0+\int_0^r\int_0^{s_*}\left(\frac{s}{s_*}\right)^{N-1}\left(-s^{\sigma}h(s)+U(s)h(s)+F(s)\right)\,ds\,ds_*\\
&=h_0+\int_0^rs^{N-1}\left(-s^{\sigma}h(s)+U(s)h(s)+F(s)\right)\int_s^{r}s_*^{1-N}\,ds_*\,ds\\
&=h_0+\int_0^rs^{N-1}\left(-s^{\sigma}h(s)+U(s)h(s)+F(s)\right)\omega_{N}(r,s)\,ds,
\end{split}
\end{equation}
with
\begin{equation}\label{interm3}
\begin{split}
&\omega_N(r,s):=\omega_N(r)-\omega_N(s), \quad N\geq1,\\
&\omega_1(s):=s, \quad \omega_2(s):=\ln\,s,\\
&\omega_N(s):=-\frac{s^{2-N}}{N-2}, \quad N\geq3.
\end{split}
\end{equation}
We split the remaining part of the proof into three cases according to the previous expressions of $\omega_N$.

\medskip

\noindent \textbf{Case 1: $N=1$.} Let $R_0\in(0,R)$ and $r\in [0,R_0]$. Taking absolute values in~\eqref{interm2}, we obtain
\begin{equation*}
	|h(r)|\leq|h_0|+\int_0^r\left((s^{\sigma}+|U(s)|)|h(s)|+ |F(s)|\right)(r-s)\,ds
\end{equation*}
and thus, recalling that $F$ is integrable on $[0,R_0]$,
\begin{equation}\label{gronw1}
	|h(r)|\leq|h_0|+ R_0\|F\|_{L^{1}(0,R_0)} + R_0 \int_0^r (s^{\sigma}+|U(s)|)|h(s)|\,ds, \quad r\in [0,R_0].
\end{equation}
Setting
\begin{equation*}
	C_0(R_0) := \blm{\max\{1,R_0\} \;\text{ and }\; C_1(R_0) := C_0(R_0) \big(|h_0| + \|F\|_{L^{1}(0,R_0)}\big),}
\end{equation*}
the estimate~\eqref{gronw1} becomes
\begin{equation*}
	|h(r)|\leq C_1(R_0) +R_0 \int_{0}^r(s^{\sigma}+|U(s)|)|h(s)|\,ds, \quad r\in [0,R_0].
\end{equation*}
Since $\sigma>-1$ and $U\in \blm{L^1(0,R_0)}$, an application of Gronwall's lemma leads us to
\begin{equation}\label{esth1}
	|h(r)|\leq C_1(R_0) \exp\left\{R_0 \left(\frac{1}{\sigma+1}r^{\sigma+1}+\|U\|_{L^1(0,r)}\right) \right\}, \quad r\in [0,R_0].
\end{equation}

\medskip

\noindent \textbf{Case 2: $N=2$.} Let $R_0\in(0,R)$ \blm{and $r\in [0,R_0]$.} Taking absolute values in~\eqref{interm2}, we obtain
\begin{equation*}
	|h(r)|\leq |h_0|+\int_0^r s\left(s^{\sigma}|h(s)|+|U(s)| |h(s)|+ |F(s)|\right)(\ln\,r-\ln\,s)\,ds, \quad r\in [0,R_0].
\end{equation*}
Fix now $\epsilon\in(0,\sigma+2)\cap (0,1)$. Since
\begin{equation*}
	0 \le s^\epsilon (\ln{r}-\ln{s}) \le s^\epsilon(\ln{R_0}-\ln{s}) \le \frac{R_0^\epsilon}{e\epsilon} \le \frac{R_0^\epsilon}{\epsilon}, \quad r\in [0,R_0],
\end{equation*}
we infer from the previous inequality and the integrability of $F$ on $[0,R_0]$ that
\begin{align}
	|h(r)| & \leq |h_0|+\int_0^r s^{1-\epsilon}\left((s^{\sigma}+|U(s)|)|h(s)|+|F(s)|\right) s^{\epsilon}(\ln\,R_0-\ln\,s)\,ds \nonumber\\
	& \le |h_0| + \frac{R_0}{\epsilon} \|F\|_{L^1(0,R_0)} + \frac{R_0^\epsilon}{\epsilon} \int_0^r \left( s^{\sigma+1-\epsilon} + s^{1-\epsilon} |U(s)| \right) |h(s)|\,ds \label{gronw2}
\end{align}
for any $r\in [0,R_0]$. Setting
\begin{equation*}
	C_0(R_0) := \frac{\max\{R_0^\epsilon , R_0\}}{\epsilon} \;\text{ and }\; C_1(R_0) := \blm{C_0(R_0) \big(|h_0| + \|F\|_{L^1(0,R_0)}\big)}
\end{equation*}
and observing that
\begin{equation*}
	s\mapsto s^{\sigma+1-\epsilon}\in L^1(0,R_0) \;\text{ and}\; s\mapsto s^{1-\epsilon} U(s)\in L^{1}(0,R_0)
\end{equation*}
due to $\sigma+1-\epsilon>-1$ and $1-\epsilon>0$, we apply Gronwall's lemma to \eqref{gronw2} and obtain
\begin{equation}\label{esth2}
	|h(r)|\leq C_1(R_0) \exp\left\{ \frac{R_0^\epsilon}{\epsilon} \left(\frac{r^{\sigma+2-\epsilon}}{\sigma+2-\epsilon}+r^{1-\epsilon}\|U\|_{L^1(0,\blm{r})}\right) \right\}, \quad r\in [0,R_0].
\end{equation}

\medskip

\noindent \textbf{Case 3: $N\geq3$.} Let $R_0\in(0,R)$ and $r\in(0,R_0)$. Taking absolute values in~\eqref{interm2}, we deduce from the integrability of $F$ on $(0,R_0)$ that
\begin{align}
	|h(r)|&\leq|h_0|+\int_0^rs^{N-1}\left(s^{\sigma}|h(s)|+|U(s)| |h(s)| + |F(s)|\right)\frac{s^{2-N}}{N-2}\,ds \nonumber \\
	&\leq |h_0| + R_0 \|F\|_{L^{1}(0,R_0)} + \int_0^r \left( s^{\sigma+1}+s|U(s)| \right) |h(s)| \,ds.\label{gronwN}
\end{align}
Setting
\begin{equation*}
	C_0(R_0) := \max\{1,R_0\} \;\text{ and }\; C_1(R_0) := \blm{C_0(R_0) \big(|h_0| + \|F\|_{L^{1}(0,R_0)}\big),}
\end{equation*}
we infer from~\eqref{gronwN} and an application of Gronwall's lemma that
\begin{equation}\label{esthN}
	|h(r)|\leq C_1(R_0) \exp\left\{ \left(\frac{r^{\sigma+2}}{\sigma+2}+r\|U\|_{L^1(0,r)}\right) \right\}, \quad r\in [0,R_0].
\end{equation}
The proof is completed by observing that~\eqref{esth1}, \eqref{esth2} and~\eqref{esthN} imply~ \eqref{esth} in the three cases, according to the definition of $C_0(R_0)$. Moreover, when~\eqref{global} holds true, we can take $R_0=R$ in the previous computations and thus obtain the bound~\eqref{esth} for $r\in[0,R]$.
\end{proof}

The following uniqueness result is an obvious consequence of Lemma~\ref{lem.uniq}.
\begin{lemma}[Uniqueness]\label{cor.uniq}
Let $h\in C^2((0,R))\cap C([0,R))$ be a solution to the Cauchy problem~\eqref{singoper} \blm{with $h_0=0$, $U\in L^{1}_{\rm loc}([0,R))$ and $F\equiv0$. Then} $h\equiv0$.
\end{lemma}

\begin{proof}
We simply observe that $F\equiv0$ and $h_0=0$ imply that $C_1(R_0)=0$ for all $R_0\in (0,R)$ and then $h\equiv0$.
\end{proof}

The next result ensures the existence of at least one solution in $C^2((0,R))\cap C([0,R))$ to the Cauchy problem~\eqref{singoper} without source term.
\begin{lemma}[Well-posedness]\label{lem.exist}
Assume that $U\in L^{1}_{\rm loc}([0,R))$ and $F\equiv0$. For any $h_0\in\mathbb{R}$, there exists a unique solution $h\in C([0,R))\cap C^2((0,R))$ to the Cauchy problem~\eqref{singoper}. Moreover, if $U\in L^{\infty}_{\rm loc}([0,R))$, then
\begin{equation}\label{exph}
	h'(r)=-\frac{h_0}{\sigma+N}r^{\sigma+1}+o(r^{\sigma+1})\blm{,} \quad {\rm as} \ r\to0.
\end{equation}
\end{lemma}

\begin{proof}
The proof is based on the Banach fixed point theorem. In order to apply it, we let $R_0\in (0,R)$ to be chosen later and consider the subset 
\begin{equation*}
	X_{R_0} := \left\{ h\in C([0,R_0])\ :\ h(0)= h_0 \;\text{ and }\; \|h\|_{L^\infty(0,R_0)} \le 2|h_0| \right\}
\end{equation*} 
of $C([0,R_0])$, which is a complete metric space for the distance induced by the $L^\infty$-norm. We next define the operator $\Lambda$ which maps any function $h\in X_{R_0}$ to
\begin{equation}\label{interm4}
	\Lambda[h](r):=h_0+\int_0^rs^{N-1}h(s)\left(-s^{\sigma}+U(s)\right)\omega_N(r,s)\,ds, \quad r\in[0,R_0],
\end{equation}
where $\omega_N$ is defined in~\eqref{interm3}. We first prove that $\Lambda[h]$ is well-defined and belongs to $C([0,R_0])$ for $h\in X_{R_0}$. As in the proof of Lemma~\ref{lem.uniq}, we have
\begin{equation}\label{interm5}
\blm{0 \le s^{\sigma+N-1}\omega_N(r,s)} \leq
\begin{cases}
   s^{\sigma}R_0, & \mbox{if } N=1,\\
   s^{\sigma+1-\epsilon}s^{\epsilon}(\ln\,R_0-\ln\,s) \le \displaystyle{\frac{R_0^\epsilon}{\epsilon} s^{\sigma+1-\epsilon}}, & \mbox{if } N=2,\\
   s^{\sigma+1}, & \mbox{if } N\geq3,
\end{cases}
\end{equation}
\blm{for $0\le s \le r \le R_0$,} where $\epsilon\in(0,\sigma+2)\cap (0,1)$ is arbitrary but fixed in dimension $N=2$. We deduce from~\eqref{interm5}, the boundedness of $h$ and the integrability of $U$ that the function
\begin{equation*}
	s\mapsto s^{N-1}h(s)\left(-s^{\sigma}+U(s)\right)\omega_N(r,s)
\end{equation*}
belongs to $L^1(0,r)$ for any $r\in(0,R_0]$ and it is actually bounded from above by a function in $L^1(0,R_0)$ which is independent of $r$. A classical consequence of Lebesgue's dominated convergence theorem is that $\Lambda[h]\in C([0,R_0])$ with $\Lambda[h](0)=h_0$. We next infer from~\eqref{interm4} that 
\begin{align*}
	|\Lambda[h](r)| & \leq |h_0| + \|h\|_{L^\infty(0,R_0)} \int_0^r s^{\sigma+N-1} \omega_N(r,s)\,ds \\
	& \quad +  \|h\|_{L^\infty(0,R_0)} \sup_{s\in [0,r]}\left\{ s^{N-1} \omega_N(r,s) \right\} \|U\|_{L^{1}(0,r)},
\end{align*}
for any $r\in[0,R_0]$, and the estimates~\eqref{interm5}, along with the integrability properties of $U$, ensure that there is $R_0\in (0,R)$ such that
\begin{equation*}
	|\Lambda[h](r)|\leq |h_0| + \frac{\|h\|_{L^\infty(0,R_0)} }{2}\leq 2|h_0|, \quad \blm{r\in [0,R_0]}.
\end{equation*}
Therefore, $\Lambda[h]\in X_{R_0}$. Similarly, picking two functions $h_1\in X_{R_0}$ and $h_2\in X_{R_0}$, we deduce from~\eqref{interm4} that, for $r\in [0,R_0]$, 
\begin{equation*}
\begin{split}
	|\Lambda[h_1](r)-\Lambda[h_2](r)| & = \left| \int_0^r s^{N-1} (h_1-h_2)(s) \left(-s^{\sigma}+U(s)\right) \omega_N(r,s)\,ds \right| \\
	& \leq \|h_1-h_2\|_{L^\infty(0,R_0)} \int_0^r s^{\sigma+N-1} \omega_N(r,s)\,ds \\
	& \quad + \|h_1-h_2\|_{L^\infty(0,R_0)}  \sup_{s\in [0,r]}\left\{ s^{N-1} \omega_N(r,s) \right\} \|U\|_{L^{1}(0,r)}\\
	&\leq \frac{1}{2} \|h_1-h_2\|_{\blm{L^\infty(0,R_0)}},
\end{split}
\end{equation*}
by the previous choice of $R_0$. An application of the Banach fixed point theorem ensures that $\Lambda$ has a unique fixed point $h\in X_{R_0}$; that is, there is a unique solution $h\in X_{R_0}$ to
\begin{equation*}
	h(r)=h_0+\int_0^rs^{N-1}h(s)\left(-s^{\sigma}+U(s)\right)\omega_N(r,s)\,ds, \quad r\in[0,R_0].
\end{equation*}
Recalling from~\eqref{interm3} that $\omega_N(r,s)=\omega_N(r)-\omega_N(s)$, we further obtain that, for $r\in[0,R_0]$,
\begin{equation*}
\begin{split}
	h(r)&=h_0+\omega_N(r)\int_0^rs^{N-1}h(s)\left(-s^{\sigma}+U(s)\right)\,ds\\
	& \quad -\int_0^rs^{N-1}\omega_N(s)h(s)\left(-s^{\sigma}+U(s)\right)\,ds,
\end{split}
\end{equation*}
both integrals being finite in view of the fact that $\sigma>\max\{-N,-2\}$, the integrability of $U$ and the definition of $\omega_N$. Consequently, $h\in C^1((0,R_0))$ with
\begin{equation*}
	h'(r)=\omega_N'(r)\int_0^r s^{N-1} h(s) \left(-s^{\sigma}+U(s)\right)\,ds, \quad r\in (0,R_0).
\end{equation*}
Since $\omega_N'(r)=r^{1-N}$ for any $N\geq1$, we find
\begin{equation}\label{interm6}
	r^{N-1}h'(r)=\int_0^rs^{N-1}h(s)\left(-s^{\sigma}+U(s)\right)\,ds, \quad r\in(0,R_0),
\end{equation}
which implies, by differentiating once more with respect to $r$, that $h$ solves~\eqref{L0} on $(0,R_0)$. Moreover, \eqref{interm6} also implies that
\begin{equation*}
	\lim\limits_{r\to0}r^{N-1}h'(r)=0,
\end{equation*}
hence $h$ satisfies the initial conditions~\eqref{ic}. We have thus shown that $h$ is a solution to the Cauchy problem~\eqref{singoper}. In addition, $h$ cannot become unbounded in $(0,R)$ according to Lemma~\ref{lem.uniq}. Classical arguments then allow us to extend $h$ to the whole interval $(0,R)$ and thereby complete the proof of the well-posedness of the Cauchy problem~\eqref{singoper} on $(0,R)$, the uniqueness being granted by Lemma~\ref{cor.uniq}.

Finally, assuming further that $U\in L^\infty_{\rm loc}([0,R))$, it follows from the continuity of $h$ on $[0,R)$ that
\begin{equation*}
	\int_0^r s^{\sigma+N-1}h(s)\,ds \sim \frac{h_0r^{\sigma+N}}{\sigma+N}, \quad \int_0^{r}s^{N-1}U(s)h(s)\,ds = O(r^N), \quad {\rm as} \ r\to 0,
\end{equation*}
\blm{and, taking into account the negativity of $\sigma$, we infer by substituting the previous expansions in~\eqref{interm6}} that
\begin{equation*}
	r^{N-1}h'(r)\sim-\frac{h_0}{\sigma+N} r^{\sigma+N}\blm{,} \quad {\rm as} \ r\to0,
\end{equation*}
which leads to~\eqref{exph}.
\end{proof}

We end this section of general results with a well-posedness result related to Lemma~\ref{lem.uniq}, which is fundamental in proving that solutions to Eq.~\eqref{eq1.rad} vanishing at a point $r_0>0$ with a non-zero derivative can be extended to \blm{a right neighborhood of $r_0$}.
\begin{proposition}\label{prop.neg}
Let $r_0>0$ and $b\in\mathbb{R}\setminus\{0\}$. The Cauchy problem 
\begin{subequations}\label{CPneg}
\begin{equation}\label{eq.neg}
	g''(r)+\frac{N-1}{r}g'(r)+(|g|^{(1-m)/m}g)(r)-r^{\sigma}g(r)=0, \quad r>r_0,
\end{equation}
\begin{equation}\label{icneg}
	g(r_0)=0, \quad g'(r_0)=b,
\end{equation}
\end{subequations}
has a unique solution $g\in C^2([r_0,R_0))$ satisfying $bg>0$ in $(r_0,R_0)$ and defined on a maximal interval $(r_0,R_0)$ for some $R_0\in(r_0,\infty]$. In addition, if $R_0<\infty$, then
\begin{equation}
	\lim_{r\to R_0} \left( |g(r)| + |g'(r)| + \frac{1}{|g(r)|} \right) = \infty. \label{uc}
\end{equation} 
\end{proposition}
\begin{proof}
Obviously, the Cauchy-Lipschitz theorem cannot be applied to~\eqref{eq.neg} due to the singularity of the function $z\mapsto |z|^{(1-m)/m}z$ at $z=0$. However, the Cauchy-Peano theorem implies that the Cauchy problem \eqref{CPneg} has at least one local solution $g\in C^2([r_0,\bar{r}_0])$ for some $\bar{r}_0\in (r_0,\infty)$.  At this point, if $g(\bar{r}_0)=0$, then we set $R_0=\bar{r}_0$ and~\eqref{uc} is satisfied. If $g(\bar{r}_0)\neq 0$, then the Cauchy-Lipschitz theorem allows us to extend $g$ to a larger interval $[r_0,R_0)$ and the finiteness of $R_0$ entails~\eqref{uc}.

As for the uniqueness, let us \blm{assume} that~\eqref{CPneg} admits two solutions $g_1\in C^2([r_0,R_1))$ and $g_2\in C^2([r_0,R_2))$ for some $R_2\ge R_1>r_0$ such that $bg_i$ is positive on $(r_0,R_i)$, $i=1,2$. On the one hand, there is $\varrho>r_0$ such that 
\begin{equation}\label{neg1}
|g_i(r)|\geq\frac{|b|}{2}(r-r_0), \quad r\in[r_0,\varrho), \quad i=1,2.
\end{equation}
On the other hand, we argue exactly as in the deduction of~\eqref{interm2} in Lemma~\ref{lem.uniq} (but with $g_2-g_1$ in place of $h$, and zero in place of $h_0$) to obtain the identity
\begin{equation}\label{neg2}
	(g_2-g_1)(r)=\int_{r_0}^rs^{N-1}\left[-s^{\sigma}(g_2-g_1)(s)+D(s)(g_2-g_1)(s)\right]\omega_N(r,s)\,ds, \quad r\in [r_0,R_1),
\end{equation}
with 
\begin{equation*}
D(s):=
\begin{cases}
	\displaystyle{\frac{(|g_1|^{(1-m)/m}g_1)(s)-(|g_2|^{(1-m)/m}g_2)(s)}{g_1(s)-g_2(s)}}, & \mbox{if } g_1(s)\neq g_2(s), \\[2mm]
	\displaystyle{\frac{1}{m}|g_1|^{(1-m)/m}(s)}, & \mbox{if } g_1(s)=g_2(s).
\end{cases}
\end{equation*}
Thanks to~\eqref{neg1}, we have
\begin{equation*}
	0\leq D(s)\leq\frac{1}{m}\left(\frac{|b|}{2}\right)^{(1-m)/m}(s-r_0)^{(1-m)/m}, \quad s\in[r_0,\varrho],
\end{equation*}
hence $D\in L^1(r_0,\varrho)$. In view of this property, we argue as in the proof of Lemma \ref{lem.uniq} starting from \eqref{neg2} to deduce that $g_1=g_2$ on $[r_0,\varrho]$. Next, since $g_1(\varrho)=g_2(\varrho)\neq 0$ and $g_1'(\varrho)=g_2'(\varrho)$, we deduce from the Cauchy-Lipschitz theorem that $g_1=g_2$ on $[\varrho,R_1)$, which in turn implies that $R_1=R_2$ according to~\eqref{uc}.
\end{proof}

\section{Well-posedness of~\eqref{CPVs}}\label{app.B}

In this section, we apply the general results designed in Appendix~\ref{app.A} in order to show that~\eqref{CPVs} is well-posed, \blm{as stated in Proposition~\ref{prop.wpv}}.
\begin{proposition}\label{prop.CP}
Let $m$ and $\sigma$ be as in~\eqref{range.exp} and $a\in(0,\infty)$. There are $S(a)\in (0,\infty]$ and a unique  solution 
\begin{equation*}
	V(\cdot;a)\in C([0,S(a))) \cap C^2((0,S(a)))
\end{equation*}
to~\eqref{CPVs} satisfying~\eqref{nd}; that is, $\big(V(r;a),V'(r;a)\big)\ne (0,0)$ for all $r\in (0,S(a))$. Moreover, if $S(a)<\infty$, then 
\begin{equation}
	\lim\limits_{r\to S(a)} V(r;a) = \lim\limits_{r\to S(a)} V'(r;a) = 0, \label{vanishing}
\end{equation}
and $V(\cdot;a)$ can be extended to a global $C^2$-smooth solution to~\eqref{CPVs} on $(0,\infty)$ by setting $V(r;a)=0$ for $r\in [S(a),\infty)$.
\end{proposition}

The proof of Proposition~\ref{prop.CP} builds upon the preliminary results established in Appendix~\ref{app.A} and involves several steps. We begin with the local well-posedness of~\eqref{CPVs} within the framework of positive solutions.
\begin{lemma}\label{lem.CPloc}
Let $m$ and $\sigma$ be as in~\eqref{range.exp} and $a\in(0,\infty)$. Then there is a unique positive solution
\begin{equation*}
	V(\cdot;a)\in C([0,Z_0(a))) \cap C^2((0,Z_0(a)))
\end{equation*}
to~\eqref{CPVs} defined on a maximal interval $(0,Z_0(a))$ with $Z_0(a)\in (0,\infty]$. Moreover, if $Z_0(a)<\infty$, then 
\begin{equation}
	\lim\limits_{r\to Z_0(a)} \left( V(r;a) + |V'(r;a)| + \frac{1}{V(r;a)}\right)=\infty. \label{alt}
\end{equation}

\end{lemma}
\begin{proof}
We first prove the \textbf{uniqueness} statement. Assume that there exist $R_0>0$ and two positive solutions $V_i\in C([0,R_0))\cap C^2((0,R_0))$ to the Cauchy problem~\eqref{CPV}, $i=1,2$. Setting $\overline{V}:=V_1-V_2$, it follows from~\eqref{ica} that $\overline{V}$ satisfies the initial conditions \eqref{ic} with $h_0=0$. Moreover, replacing $V_1$ and $V_2$ in~\eqref{eq1.rad} and subtracting the resulting equations, we find
\begin{equation}\label{interm22}
	\overline{V}''(r)+\frac{N-1}{r}\overline{V}'(r)+r^{\sigma}\overline{V}(r)-D(r)\overline{V}(r)=0, \quad r\in(0,R_0),
\end{equation}
where
\begin{equation*}
	D(r):=\begin{cases}
        \displaystyle{\frac{V_1^{1/m}(r)-V_2^{1/m}(r)}{V_1(r)-V_2(r)}}, & \mbox{if } V_1(r)\neq V_2(r), \\[4mm]
        \displaystyle{\frac{1}{m}V_1^{(1-m)/m}(r)}, & \mbox{if } V_1(r)=V_2(r).
    \end{cases}
\end{equation*}
Since
\begin{equation*}
	0 \le D(r) \le \frac{1}{m} \min\{V_1(r),V_2(r)\}^{(1-m)/m}, \quad r\in (0,R_0),
\end{equation*}	
the positivity of $V_1$ and $V_2$ warrants that $D\in L_{\mathrm{loc}}^\infty([0,R_0))$. We then apply Lemma~\ref{cor.uniq} to~\eqref{interm22}, with $h_0=0$, $U=D$ and $F=0$, to deduce that $\overline{V}\equiv0$ on $[0,R_0)$ and the uniqueness is established.

\medskip

We now turn to the \textbf{existence} statement. We proceed as in the proof of Lemma~\ref{lem.exist} and introduce  the subset
\begin{equation*}
	E_{R_0} := \left\{ h\in C([0,R_0])\ :\ h(0)= \blm{a} \;\text{ and }\; 0 \le h \le 2a \right\}
\end{equation*} 
of $C([0,R_0])$, which is a complete metric space for the distance induced by the $L^\infty$-norm, and the operator $\Lambda$ defined on $E_{R_0}$ by
\begin{equation*}
	\Lambda[V](r) := a + \int_0^r s^{N-1}V(s) \left(-s^{\sigma} + V^{(1-m)/m}(s) \right) \omega_N(r,s)\,ds, \quad r\in [0,R_0],
\end{equation*}	
where $\omega_N(s,r)$ is defined in \eqref{interm3}. \blm{We first prove a non-negativity-preserving property of $\Lambda$ for $R_0>0$ sufficiently small. Indeed, introducing 
\begin{equation*}
	c_1:=\frac{1}{\sigma+1}, \quad c_2:=\frac{1}{\epsilon(\sigma+2-\epsilon)}, \quad c_N:=\frac{1}{\sigma+2}, \quad {\rm for} \ N\geq3,
\end{equation*}
for some fixed $\epsilon\in(0,\sigma+2)\cap(0,1)$, we choose $R_0\le (2c_N)^{-1/(\sigma+2)}$, so that~\eqref{interm5} ensures that
\begin{equation*}
	\Lambda[V](r) \ge a - 2a \int_0^r s^{\sigma+N-1} \omega_N(r,s)\,ds \ge a \left( 1 - 2 c_N R_0^{\sigma+2} \right) \ge 0 
\end{equation*}
for $r\in [0,R_0]$ and $V\in E_{R_0}$. We then argue as in the proof of Lemma~\ref{lem.exist} to show that $\Lambda$ is a contraction on $E_{R_0}$ with respect to the $L^{\infty}$-norm for $R_0>0$ sufficiently small.} This fixed point argument proves the existence of $V(\cdot;a)$ on an interval $[0,R_0]$. We next infer from \cite[Theorem~7.6]{Amann} that the solution $V(\cdot;a)$ can be extended to a maximal interval $(0,Z_0(a))$ with $Z_0(a)\in (0,\infty]$ such that the alternative~\eqref{alt} holds true.
\end{proof}

We next discard the possibility of $V(\cdot;a)$ or $V'(\cdot;a)$ to become unbounded as $r\to Z_0(a)$ when $Z_0(a)<\infty$. To this end, we first derive an energy inequality for~\eqref{eq1.rads} in the spirit of Lemma~\ref{lem.energy}.
\begin{lemma}\label{lem.energyb}
Let $W\in C^1([r_1,r_2))\cap C^2((r_1,r_2))$ be a solution to~\eqref{eq1.rads} on $(r_1,r_2)$ for some $0<r_1<r_2$. Introducing the energy
\begin{subequations}\label{energyb2}
\begin{equation}\label{energyb2a}
	E_W(r):=\frac{1}{2}(W')^2(r)+\frac{1}{2}r^{\sigma}W^2(r)-\frac{m}{m+1}|W|^{(m+1)/m}(r), \quad r\in [r_1,r_2),
\end{equation}
there holds
\begin{equation}
	E_W(r) \le E_W(r_1), \quad r\in [r_1,r_2), \label{energyb2b}
\end{equation}
and 
\begin{equation}
	(W')^2(r)+\frac{r_2^{\sigma}}{2}W^2(r) \le 2 E_W(r_1) + \big(2 r_2^{-\sigma}\big)^{(m+1)/(m-1)}, \quad r\in [r_1,r_2). \label{energyb2c}
\end{equation}
\end{subequations}
\end{lemma}
\begin{proof}
It readily follows from~\eqref{eq1.rads}, \eqref{energyb2a} and the negativity of $\sigma$ that
\begin{equation}
	E_W'(r) = - \frac{N-1}{r} (W')^2(r) + \frac{\sigma r^{\sigma-1}}{2}W^2(r) \leq 0, \quad r\in [r_1,r_2), \label{energyb3}
\end{equation}
\blm{which leads to}~\eqref{energyb2b}. Moreover, since $(m+1)/m<2$, we infer from Young's inequality and the negativity of $\sigma$ that, for $r\in [r_1,r_2)$,
\begin{align*}
	E_W(r) & \ge \frac{1}{2}(W')^2(r)+\frac{r_2^{\sigma}}{2}W^2(r) - \frac{m}{m+1} \left( \frac{r_2^\sigma W^2(r)}{2} \right)^{(m+1)/2m} \left( 2 r_2^{-\sigma} \right)^{(m+1)/2m} \\
	& \ge \frac{1}{2} \left[ (W')^2(r) + \frac{r_2^{\sigma}}{2}W^2(r) - \left( 2 r_2^{-\sigma} \right)^{(m+1)/(m-1)} \right],
\end{align*}
which we combine with~\eqref{energyb2b} to obtain~\eqref{energyb2c}.
\end{proof}

Owing to Lemma~\ref{lem.energyb}, we may refine the outcome of Lemma~\ref{lem.CPloc}.
\begin{lemma}\label{lem.CPloc0}
Let $a\in (0,\infty)$ be such that $Z_0(a)<\infty$, recalling that $Z_0(a)$ is defined in Lemma~\ref{lem.CPloc}. Then
\begin{equation*}
	\lim_{r\to Z_0(a)} V(r;a) = 0 \;\;\text{ and }\;\; \lim_{r\to Z_0(a)} V'(r;a) \in\mathbb{R}.
\end{equation*}
\end{lemma}

\begin{proof}
Pick any $r_1\in (0,Z_0(a))$. For $r\in [r_1,Z_0(a))$, we infer from~\eqref{energyb2c} (with $W=V(\cdot;a)$ and $r_2=Z_0(a)$) that
\begin{equation*}
	(V')^2(r;a) + \frac{Z_0(a)^{\sigma}}{2} V^2(r;a) \le 2 E_{V(\cdot;a)}(r_1) + \big(2 Z_0(a)^{-\sigma}\big)^{(m+1)/(m-1)},
\end{equation*}
an inequality which prevents both $V(\cdot;a)$ and $|V'(\cdot;a)|$ to be unbounded as $r\to Z_0(a)$. Therefore, in view of~\eqref{alt}, we necessarily have
\begin{equation}
	\lim_{r\to Z_0(a)} V(r;a) = 0. \label{energyb5}
\end{equation}
Moreover, according to~\eqref{CPVs} and the positivity of $V(\cdot;a)$ on $(0,Z_0(a))$,
\begin{equation*}
	r^{N-1} V'(r;a) = \int_0^r s^{N-1} \left( - s^\sigma V(s;a) + V^{1/m}(s;a) \right)\,ds, \quad r\in (0,Z_0(a)).
\end{equation*}
Since $\sigma+N>0$ and $V(\cdot;a)$ is bounded on $(0,Z_0(a))$ by~\eqref{energyb5}, the above identity implies that $V'(\cdot;a)$ has a finite limit as $r\to Z_0(a)$, and the proof is complete.
\end{proof}

\begin{proof}[Proof of Proposition~\ref{prop.CP}]
Let $a\in (0,\infty)$ and recall that we have already constructed a unique positive solution $V(\cdot;a)$ to~\eqref{CPVs} on $(0,Z_0(a))$ in Lemma~\ref{lem.CPloc}. Our main purpose is now to extend this solution beyond $Z_0(a)$ in a unique way and we employ a iterative process to this end.

\medskip

\noindent\textbf{Step~1.} To begin with, two cases are to be considered, depending on whether $Z_0(a)$ is finite.

\smallskip

\noindent$\triangleright$ \textbf{Case~1:} Either $Z_0(a)=\infty$. In that case, we set $S(a)=\infty$  and $V(\cdot;a)$ is a positive solution to~\eqref{CPVs} on $(0,\infty)$.

\smallskip

\noindent$\triangleright$ \textbf{Case~2:} Or $Z_0(a) <\infty$ and Lemma~\ref{lem.CPloc0} ensures that $V(\cdot,a)$ has a $C^2$-smooth extension to $(0,Z_0(a)]$ with 
\begin{equation*}
	V(Z_0(a);a)=0, \quad V'(Z_0(a);a) = \lim_{r\to Z_0(a)} V'(r;a)\in\mathbb{R}, \quad V''(Z_0(a);a) = - \frac{N-1}{Z_0(a)} V'(Z_0(a);a).
\end{equation*}
The subsequent analysis then depends on the value of $V'(Z_0(a);a)$. 

\smallskip

$\triangleright$ \textbf{Case~2.1:} If $V'(Z_0(a);a)=0$, then we set $S(a)=Z_0(a)$ and realize that $V''(Z_0(a);a)=0$. Therefore, setting $V(r;a)=0$ for $r>Z_0(a)$, the function $V(\cdot;a)$ thus constructed is clearly a $C^2$-smooth solution to~\eqref{CPVs} on $(0,\infty)$, which is positive on $(0,Z_0(a))$ and vanishes identically on $(Z_0(a),\infty)$. 

\smallskip

$\triangleright$ \textbf{Case~2.2:} If $V'(Z_0(a);a)\ne 0$, then Proposition~\ref{prop.neg} guarantees that~\eqref{eq1.rads} has a unique solution $V_1\in C^2([Z_0(a),Z_1(a))$ satisfying $V_1(Z_0(a))=0$ and $V_1'(Z_0(a))=V'(Z_0(a);a)$, which does not vanish on $(Z_0(a),Z_1(a))$ and satisfies, either $Z_1(a)=\infty$, or $Z_1(a)<\infty$ and
\begin{equation}
	\lim_{r\to Z_1(a)} \left( |V_1(r)| + |V_1'(r)| + \frac{1}{|V_1(r)|} \right) = \infty. \label{energyb8}
\end{equation}
Using Lemma~\ref{lem.energyb} (with $W=V_1$, $r_1=Z_0(a)$ and $r_2=Z_1(a)$), along with~\eqref{energyb8}, allows us to conclude that, if $Z_1(a)<\infty$, then
\begin{equation}
	\lim_{r\to Z_1(a)} V_1(r) = 0. \label{energyb9}
\end{equation}
In addition, in that case, 
\begin{equation*}
	r^{N-1} V_1'(r) = Z_0(a)^{N-1} V'(Z_0(a);a) + \int_{Z_0(a)}^r s^{N-1} \left( - s^\sigma V_1(s) + \big(|V_1|^{(1-m)/m}V_1\big)(s) \right)\,ds,
\end{equation*}
for $r\in [Z_0(a),Z_1(a))$ and the boundedness of $V_1$ on $(Z_0(a),Z_1(a))$ implies that $V_1'(r)$ has a finite limit as $r\to Z_1(a)$. Therefore, setting $V(r;a)=V_1(r)$ for $r\in(Z_0(a),\infty)$ when $Z_1(a)=\infty$ and $V(r;a)=V_1(r)$ for $r\in(Z_0(a),Z_1(a)]$ when $Z_1(a)<\infty$ provides a $C^2$-smooth solution $V(\cdot;a)$ to~\eqref{CPVs} on $(0,Z_1(a))$ and $(0,Z_1(a)]$, respectively, and its uniqueness stems from Proposition~\ref{prop.neg} and Lemma~\ref{lem.CPloc}. 

\medskip

\noindent\textbf{Step~2.} Iterating this procedure, we find a subset $J\subset \mathbb{N}$ of non-negative integers containing at least $\{0\}$, an increasing sequence $(Z_j(a))_{j\in J}$ of real numbers and a function
\begin{equation*}
	V(\cdot;a) \in C([0,S(a))\cap C^2((0,S(a))) \;\;\text{ with }\;\; S(a) := \sup_{j\in J} Z_j(a) \in [R(a),\infty],
\end{equation*}
which is a solution to~\eqref{CPVs} on $(0,S(a))$ such that
\begin{subequations}\label{V}
\begin{equation}
	V(r;a)\ne 0 \;\text{ for all }\; r\in [0,S(a))\setminus\{Z_j(a)\ :\ j\in J\}, \label{V1}
\end{equation}
\begin{equation}
	V(\cdot;a) \;\text{ is the unique positive solution to~\eqref{CPVs} on }\; (0,Z_0(a)), \label{V2} 
\end{equation}
and, for each $j\in J$ satisfying $Z_j(a)<S(a)$, 
\begin{equation}
\begin{split}
	& V(\cdot;a) \;\text{ is the unique solution to~\eqref{eq1.rads} on }\; (Z_j(a),Z_{j+1}(a)) \\
	& \text{ satisfying }\; V(Z_j(a);a) = 0 \;\text{ and }\; V'(Z_j(a);a)\ne 0.
\end{split} \label{V3}
\end{equation}
\end{subequations}

\medskip

\noindent\textbf{Step~3.} At this point, we again face different situations, which we study separately.

\smallskip

\noindent$\triangleright$ \textbf{Case~1:} If $J$ is finite, then there is $j_a\ge 1$ such that $J=\{0,\ldots,j_a\}$ and $S(a)=Z_{j_a}(a)$. 

\smallskip

$\triangleright$ \textbf{Case~1.1:} Either $Z_{j_a}(a)=\infty$ and $V(\cdot;a)$ is a solution to~\eqref{CPVs} on $(0,\infty)$.

\smallskip

$\triangleright$ \textbf{Case~1.2:} Or $Z_{j_a}(a)<\infty$ and $V(Z_{j_a}(a);a)=V'(Z_{j_a}(a);a)=0$. We then set $V(r;a)=0$ for $r\in (Z_{j_a}(a),\infty)$ to obtain a $C^2$-smooth solution to~\eqref{CPVs} on $(0,\infty)$, which vanishes identically on $[Z_{j_a}(a),\infty)$.

\smallskip

\noindent$\triangleright$ \textbf{Case~2:} If $J=\mathbb{N}$, then we again split the analysis according to the finiteness of $S(a)$.

\smallskip

$\triangleright$ \textbf{Case~2.1:} Either $S(a)=\infty$ and $V(\cdot;a)$ is a $C^2$-smooth solution to~\eqref{CPVs} on $(0,\infty)$ (with countably many sign changes).

\smallskip

$\triangleright$ \textbf{Case~2.2:} Or $S(a)<\infty$ and 
\begin{equation}
	\lim_{j\to\infty} Z_j(a) = \sup_{j\in\mathbb{N}} Z_j(a) = S(a). \label{bV1}
\end{equation}
Owing to the finiteness of $S(a)$, Lemma~\ref{lem.energyb} (with $W=V(\cdot;a)$, $r_1=Z_0(a)$ and $r_2=S(a)$) entails that
\begin{equation*}
	(V')^2(r;a) + \frac{S(a)^\sigma}{2} V^2(r;a) \le E_{V(\cdot;a)}(Z_0(a);a) + \big(2 S(a)^{-\sigma}\big)^{(m+1)/(m-1)}
\end{equation*}
for $r\in [Z_0(a),S(a))$, whence there is $L>0$ (depending on $a$) such that
\begin{equation}
	|V(r;a)| + |V'(r;a)| \le L, \quad r\in [Z_0(a),S(a)). \label{bV2}
\end{equation}
In particular, for $j\in\mathbb{N}$ and $r\in (Z_j(a);Z_{j+1}(a))$, we infer from~\eqref{V3} and~\eqref{bV2} that
\begin{equation*}
	|V(r;a)| = |V(r;a)-V(Z_j(a);a)| \le L |r-Z_j(a)| \le L |Z_{j+1}(a) - Z_j(a)|,
\end{equation*}
which we combine with~\eqref{bV1} to conclude that
\begin{equation}
	\lim_{r\to S(a)} V(r;a) = 0. \label{bV3}
\end{equation}
Next, by Rolle's theorem, there is $\eta_j(a)\in (Z_j(a),Z_{j+1}(a))$ such that $V'(\eta_j(a);a)=0$ and it follows from~\eqref{eq1.rads} and~\eqref{bV2} that, for $r\in (\eta_j(a),\eta_{j+1}(a))$,
\begin{align*}
	r^{N-1} |V'(r;a)| & = \left| r^{N-1} V'(r;a) - \eta_j(a)^{N-1} V'(\eta_j(a);a) \right| \\
	& = \left| \int_{\eta_j(a)}^r s^{N-1} \left( - s^\sigma V(s;a) + \big(|V|^{(1-m)/m} V\big)(s;a) \right)\,ds \right| \\
	& \le \int_{\eta_j(a)}^r s^{N-1} \left( Z_0(a)^\sigma L + L^{1/m} \right)\,ds \le \frac{Z_0(a)^\sigma L + L^{1/m}}{N} \left( \eta_{j+1}(a)-\eta_j(a) \right) \\
	& \le \frac{Z_0(a)^\sigma L + L^{1/m}}{N} \left( Z_{j+2}(a)-Z_j(a) \right),
\end{align*}
hence
\begin{equation}
	\lim_{r\to S(a)} V'(r;a) = 0, \label{bV4}
\end{equation} 
after using again~\eqref{bV1}. Owing to~\eqref{bV3} and~\eqref{bV4}, we may set $V(r;a)=0$ for $r\in [S(a),\infty)$ to obtain a $C^2$-smooth solution to~\eqref{CPVs} on $(0,\infty)$, which vanishes identically on $[S(a),\infty)$.

\medskip

Gathering the outcome of the above analysis, we have proved the existence of $S(a)\in (0,\infty]$ and a unique solution $V(\cdot;a)\in C([0,S(a)))\cap C^2((0,S(a)))$ to~\eqref{CPVs} on $(0,S(a))$ which satisfies~\eqref{nd}, more precise information being actually available according to~\eqref{V}, and which can be extended to a global $C^2$-smooth solution to~\eqref{CPVs} by setting $V(r;a)=0$ for $r\in [S(a),\infty)$ when $S(a)<\infty$. 
\end{proof}

\section{Continuous dependence for~\eqref{CPVs}}\label{app.C}

In this section, we investigate the stability of $V(\cdot;a)$ with respect to $a$. We emphasize once more that the singularity \blm{of~\eqref{eq1s.rad} as} $r\to 0$ or as $V\to 0$ does not allow us to use classical continuous dependence results. The aim of this section is then to establish such results. To begin with, we gather some properties of $V(\cdot;a)$ and a quantitative estimate for $R(a)$ that will be useful later.
\begin{lemma}\label{lem.propV}
For $a\in(0,\infty)$, 
\begin{equation}\label{est.zero}
	-\frac{ar^{\sigma+1}}{\sigma+N}\leq V'(r;a)<0, \quad a-\frac{ar^{\sigma+2}}{(\sigma+2)(\sigma+N)}\leq V(r;a)\leq a, \quad r\in(0,R(a)).
\end{equation}
Moreover, for each $a\in (0,\infty)$, there is a unique positive solution $Y(a)$ to
\begin{equation}\label{unif.est}
	Y(a)^{-\sigma m/(m-1)}+\frac{a Y(a)^{\sigma+2}}{(\sigma+2)(\sigma+N)} =a,
\end{equation}
and $R(a)\geq Y(a)$. In fact, $Y\in C((0,\infty))$ and satisfies
\begin{equation}\label{propy}
\begin{split}
	& 0 < Y(a) < Y_\infty := [(\sigma+2)(\sigma+N)]^{1/(\sigma+2)}, \quad \lim_{a\to\infty} Y(a) = Y_\infty, \\
	& \lim_{a\to 0} a^{(m-1)/(m\sigma)} Y(a) = 1.
\end{split}
\end{equation}
\end{lemma}
\begin{proof}
By the definition of $R(a)$, we have $0<V(r;a)<a$ and $V'(r;a)<0$ for $r\in(0,R(a))$. We then infer from~\eqref{main.id} that
\begin{equation*}
	\frac{d}{dr}(r^{N-1}V'(r;a))\geq-r^{\sigma+N-1}V(r;a)\geq-ar^{\sigma+N-1}, \quad r\in(0,R(a)),
\end{equation*}
\blm{and, taking into account the initial condition~\eqref{ica}, an integration over $(0,r)$ gives}
\begin{equation*}
	r^{N-1}V'(r;a)\geq-\frac{ar^{\sigma+N}}{\sigma+N}, \quad r\in(0,R(a)),
\end{equation*}
which is equivalent to the first estimate in~\eqref{est.zero}. Integrating once more over $(0,r)$ gives the second estimate in~\eqref{est.zero}. 

Next, for $a\in (0,\infty)$, we introduce the function 
\begin{equation*}
	\xi_a(z) := z^{-m\sigma/(m-1)} + a \left( \frac{z}{Y_\infty} \right)^{\sigma+2}, \quad z\in (0,\infty),
\end{equation*}
and deduce from the positivity of $-\sigma$ and $\sigma+2$ that $\xi_a$ is an increasing function from $(0,\infty)$ onto $(0,\infty)$. Consequently, there is a unique $Y(a)\in (0,\infty)$ satisfying $\xi_a(Y(a))=a$ or, equivalently, the identity~\eqref{unif.est}. Now, the monotonicity of $\xi_a$, along with the explicit formula $\xi_a(Y_\infty) = a + Y_\infty^{-m\sigma/(m-1)}>a$, entails that $Y(a)<Y_\infty$. Furthermore, the \blm{identity~\eqref{unif.est}} readily ensures that $Y\in C((0,\infty))$ and, since
\begin{equation*}
	\left( \frac{Y(a)}{Y_\infty} \right)^{\sigma+2} = 1 - \frac{Y(a)^{-m\sigma/(m-1)}}{a}, \quad a\in (0,\infty),
\end{equation*}
we let $a\to\infty$ in the above identity and deduce from the just established boundedness of $Y$ that $Y(a)\to Y_\infty$ as $a\to\infty$, thereby completing the proof of the first two statements in~\eqref{propy}. As for the behavior of $Y(a)$ as $a\to 0$, we observe on the one hand that
\begin{equation*}
	\xi_a\big(a^{-(m-1)/(m\sigma)}\big) = a + \frac{a^{-(2m-\sigma-2)/(m\sigma)}}{Y_\infty^{\sigma+2}} \ge a, \quad a\in (0,\infty),
\end{equation*}
so that
\begin{equation}
	a^{-(m-1)/(m\sigma)} \ge Y(a), \quad a\in (0,\infty). \label{ya0u}
\end{equation}
On the other hand, for $\theta\in (0,1)$ and
\begin{equation*}
	0 < a < a_\theta := \left( \frac{Y_\infty}{\theta} \right)^{-m\sigma/(m-1)} \left( 1 - \theta^{-m\sigma/(m-1)} \right)^{-m\sigma/[(m-1)(\sigma+2)]},
\end{equation*}
one has
\begin{align*}
	\xi_a\big(\theta a^{-(m-1)/(m\sigma)}\big) & = a \theta^{-m\sigma/(m-1)} + a \left( \frac{\theta}{Y_\infty} \right)^{\sigma+2} a^{-(m-1)(\sigma+2)/(m\sigma)} \\
	& <a \left[ \theta^{-m\sigma/(m-1)} + \left( \frac{\theta}{Y_\infty} \right)^{\sigma+2} a_\theta^{-(m-1)(\sigma+2)/(m\sigma)} \right] \\
	& = a \left[ \theta^{-m\sigma/(m-1)} + 1 - \theta^{-m\sigma/(m-1)} \right] = a,
\end{align*}
which implies that
\begin{equation}
	\theta a^{-(m-1)/(m\sigma)} < Y(a), \quad a\in (0,a_\theta). \label{ya0l}
\end{equation}
Combining~\eqref{ya0u} and~\eqref{ya0l} completes the proof of~\eqref{propy}.

We are left to prove that $R(a)\geq Y(a)$. To this end, we observe that the function
\begin{equation*}
	\eta: r\mapsto r^{-\sigma}-V^{(m-1)/m}(r;a)
\end{equation*}
is increasing on $(0,R(a))$, with a negative value $\eta(0)=-a^{(m-1)/m}$ at $r=0$. Assume for contradiction that $\eta< 0$ on $(0,R(a))$. Then $V(R(a);a)\ge R(a)^{-\sigma m/(m-1)}>0$ and we infer from~\eqref{ica} and~\eqref{main.id} that 
\begin{equation*}
	R(a)^{N-1}V'(R(a);a) = \int_0^{R(a)} s^{\sigma+N-1} \eta(s) V^{1/m}(\blm{s;a})\,ds < 0,
\end{equation*}
which contradicts the definition~\eqref{Ra} of $R(a)$. Therefore, there is a unique $z(a)\in (0,R(a))$ such that
\begin{equation*}
\eta(r) = r^{-\sigma}-V^{(m-1)/m}(r;a)\begin{cases}
                              <0, & \mbox{if } r\in(0,z(a)), \\
                              >0, & \mbox{if } r\in(z(a),R(a)).
                            \end{cases}
\end{equation*}
Combining~\eqref{est.zero} and the definition of $z(a)$ gives
\begin{equation*}
	a-\frac{az(a)^{\sigma+2}}{(\sigma+2)(\sigma+N)}\leq V(z(a);a)=z(a)^{-m\sigma/(m-1)},
\end{equation*}
or, equivalently, 
\begin{equation*}
	\xi_a(z(a)) \geq a = \xi_a(Y(a)).
\end{equation*}
The monotonicity of $\blm{\xi}_a$ then implies that $z(a)\geq Y(a)$ and thus $R(a)\geq Y(a)$, completing the proof.
\end{proof}

A direct consequence of Lemma~\ref{lem.propV} is a negative upper bound on $V'(\cdot;a)$ in a neighborhood of $r=0$, which depends continuously on $a\in (0,\infty)$.
\begin{corollary}\label{cor.propV}
Let $a\in (0,\infty)$. Then,
\begin{equation*}
	V'(r;a) \le -\frac{r^{\sigma+1}}{2(\sigma+N)} Y(a)^{-m\sigma/(m-1)}, \quad r\in \left( 0, 2^{1/\sigma}Y(a)\right).
\end{equation*}
\end{corollary}
\begin{proof}
Let $r\in \left( 0,2^{1/\sigma}Y(a)\right)\subset (0,Y(a))$. Since $\sigma+2>0$, we infer from~\eqref{est.zero} and~\eqref{unif.est} that
\begin{equation*}
	V(s;a) \ge a - \frac{a Y(a)^{\sigma+2}}{(\sigma+2)(\sigma+N)} = Y(a)^{-\sigma m/(m-1)}, \quad s\in [0,r],
\end{equation*}
from which, \blm{using as well the positivity of $V(\cdot;a)$ on $(0,Y(a))$,} we deduce that
\begin{align*}
	r^{N-1} V'(r;a) & = \int_0^r s^{\sigma+N-1} V(s;a) \left[ - 1 + s^{-\sigma} V^{(1-m)/m}(s;a) \right]\,ds \\
	& \le \int_0^r s^{\sigma+N-1} V(s;a) \left[ - 1 + s^{-\sigma} Y(a)^\sigma \right]\,ds \\
	& \le - \frac{1}{2} \int_0^r s^{\sigma+N-1} V(s;a)\,ds \le -\frac{r^{\sigma+N}}{2(\sigma+N)} Y(a)^{-m\sigma/(m-1)},
\end{align*}
as claimed.
\end{proof}

After this preparation, we are in a position to show the continuous dependence of $V(\cdot;a)$ in a neighborhood of $r=0$.
\begin{lemma}\label{lem.loccd}
Fix $a_{\min}>0$. There are $Y_{\min}>0$ and $A_{\min}>0$ depending only on $a_{\min}$ such that
\begin{equation}\label{contg}
\begin{split}
	\sup\limits_{r\in[0,Y_{\min}]}|V(r;a)-V(r;a_0)| & \le A_{\min} |a-a_0|, \quad (a_0,a)\in (a_{\min},\infty)^2, \\
	\sup\limits_{r\in[0,Y_{\min}]} \left\{ r^{N-1} |V'(r;a)-V'(r;a_0)| \right\} & \le A_{\min} |a-a_0|, \quad (a_0,a)\in (a_{\min},\infty)^2.
\end{split}
\end{equation}
\end{lemma}

\begin{proof}
According to Lemma~\ref{lem.propV} , there is $Y_{\min}\in (0,Y_\infty)$ such that
\begin{equation*}
	Y_{\min}\leq Y(a) \le R(a), \quad a\in [a_{\min},\infty),
\end{equation*}
and
\begin{equation*} 													
	a-\frac{a Y_{\min}^{\sigma+2}}{Y_\infty^{\sigma+2}} = a \left[ 1 - \left( \frac{Y_{\min}}{Y_\infty} \right)^{\sigma+2} \right] \geq \delta_{\min} := a_{\min} \left[ 1 - \left( \frac{Y_{\min}}{Y_\infty} \right)^{\sigma+2} \right] >0, \quad a\in [a_{\min},\infty).
\end{equation*}
In particular, $V(\cdot;a)$ is well-defined on $[0,Y_{\min}]$ for all $a\in [a_{\min},\infty)$ and it follows from~\eqref{est.zero} that
\begin{equation}\label{cont1}
	V(r;a)\ge a-\frac{a Y_{\min}^{\sigma+2}}{Y_\infty^{\sigma+2}} \ge \delta_{\min}>0, \quad (r,a)\in(0,Y_{\min})\times [a_{\min},\infty).
\end{equation}
Consider now $(a_0,a)\in (a_{\min},\infty)^2$. Setting $W(r):=V(r;a)-V(r;a_0)$ for $r\in [0,Y_{\min}]$, we realize that $W$ solves the initial value problem
\begin{subequations}\label{CPW}
\begin{equation}\label{cont2}
	W''(r)+\frac{N-1}{r}W'(r)+r^{\sigma}W(r)-D(r;a,a_0)W(r)=0, \quad r\in(0,Y_{\min}),
\end{equation}
\begin{equation}\label{icw}
	W(0)=a-a_0, \quad \lim\limits_{r\to 0} r^{N-1}W(r)=0,
\end{equation}
\end{subequations}
with
\begin{equation}
D(r;a,a_0):=\begin{cases}
	\displaystyle{\frac{V^{1/m}(r;a)-V^{1/m}(r;a_0)}{V(r;a)-V(r;a_0)},} & \mbox{if } V(r;a)\neq V(r;a_0), \\[4mm]
	\displaystyle{\frac{1}{m}V^{(1-m)/m}(r;a_0),} & \mbox{if } V(r;a)=V(r;a_0).
\end{cases} \label{D}
\end{equation}
For $r\in [0,Y_{\min}]$, we infer from~\eqref{cont1} that
\begin{equation*}
\begin{split}
	\big|V^{1/m}(r;a)-V^{1/m}(r;a_0)\big| &=\frac{1}{m}\left|\int_{V(r;a_0)}^{V(r;a)}z^{(1-m)/m}\,dz\right|\\
	&\leq\frac{1}{m}|V(r;a)-V(r;a_0)|\min\{V(r;a),V(r;a_0)\}^{(1-m)/m} \\
	&\leq\frac{\delta_{\min}^{(1-m)/m}}{m}|V(r;a)-V(r;a_0)|
\end{split}
\end{equation*}
and conclude that
\begin{equation}\label{cont3}
	0\leq D(r;a,a_0)\leq\frac{\delta_{\min}^{(1-m)/m}}{m} \le \delta_{\min}^{(1-m)/m}, \quad r\in(0,Y_{\min}).
\end{equation}
We infer from Lemma~\ref{lem.uniq} (with $U=D(\cdot;a,a_0)$ and $F=0$) applied to the Cauchy problem~\eqref{CPW} on $(0,Y_{\min})$ and~\eqref{cont3} that
\begin{align*}
|W(r)| & \leq |a-a_0| \exp\left\{ C_0(Y_{\min}) \left(\frac{Y_{\min}^{\sigma_N}}{\sigma_N} + Y_{\min}^{\sigma_N-\sigma-1} \|D(\cdot;a,a_0)\|_{L^1(0,Y_{\min})} \right) \right\} \\
& \leq |a-a_0| \exp\left\{ C_0(Y_{\min}) \left(\frac{Y_{\min}^{\sigma_N}}{\sigma_N} + \delta_{\min}^{(1-m)/m} Y_{\min}^{\sigma_N-\sigma} \right) \right\}, 
\end{align*}
which gives
\begin{equation}\label{cont4}
	\sup\limits_{r\in[0,Y_{\min}]}|V(r;a)-V(r;a_0)| \le A |a-a_0|, \quad (a_0,a)\in (a_{\min},\infty)^2,
\end{equation}
with 
\begin{equation*}
	A := \exp\left\{ C_0(Y_{\min}) \left(\frac{Y_{\min}^{\sigma_N}}{\sigma_N} + \delta_{\min}^{(1-m)/m} Y_{\min}^{\sigma_N-\sigma} \right) \right\}.
\end{equation*}
We next infer from~\eqref{CPW}, \eqref{cont3} and~\eqref{cont4} that, for $r\in (0,Y_{\min}]$, 
\begin{align*}
	\big|r^{N-1} W'(r)\big| & = \left| \int_0^r s^{N-1} \big( D(s;a,a_0) - s^\sigma \big) W(s)\,ds \right| \\
	& \le A |a-a_0| \int_0^r \left( \delta_{\min}^{(1-m)/m} s^{N-1} + s^{\sigma+N-1}\right)\,ds \\
	& \le A |a-a_0| \left( \delta_{\min}^{(1-m)/m} \frac{Y_{\min}^N}{N} + \frac{Y_{\min}^{\sigma+N}}{\sigma+N} \right),
\end{align*}
from which the continuous dependence of $r\mapsto r^{N-1} V'(r;a)$ with respect to $a\in (a_{\min},\infty)$ in a neighborhood of $r=0$ which only depends on $a_{\min}$ follows. Recalling~\eqref{cont4}, we have proved~\eqref{contg}.
\end{proof}

To extend the continuous dependence with respect to $a$ to an (almost) arbitrary interval, we derive estimates on $V(\cdot;a)$, which depend mildly on $a>0$.
\begin{lemma}\label{lem.gb}
Let $0<a_{\min}<a_{\max}$ and \blm{$\rho>Y_{\min}$, with $Y_{\min}>0$ defined in Lemma~\ref{lem.loccd}.} There is $L(a_{\min},a_{\max},\rho)>0$ depending only on $a_{\min}$, $a_{\max}$ and $\rho$ such that, for $a\in [a_{\min},a_{\max}]$,
\begin{equation}\label{contgb}
	|V(r;a)|  + \blm{|V'(r;a)| +|V''(r;a)|} \le L(a_{\min},a_{\max},\rho), \quad r\in \blm{[Y_{\min},\rho].}
\end{equation}
\end{lemma}

\begin{proof}
Let $a\in [a_{\min},a_{\max}]$. Since $V(\cdot;a)$ solves~\eqref{eq1.rads} on $(Y_{\min},\rho)$, \blm{we} infer from~\eqref{energyb2c} that, for $r\in [Y_{\min},\rho]$, 
\begin{align*}
	(V'(r;a))^2 & + \frac{\rho^\sigma}{2} V^2(r;a) \le 2 E_{V(\cdot;a)}(Y_{\min}) + \big(2\rho^{-\sigma}\big)^{(m+1)/(m-1)} \\
	& \le (V'(Y_{\min};a))^2 + Y_{\min}^\sigma V^2(Y_{\min};a) + \big(2\rho^{-\sigma}\big)^{(m+1)/(m-1)} \\
	& \le 2 \big|V'(Y_{\min};a)-V'(Y_{\min};a_{\min})\big|^2 + 2 \big|V'(Y_{\min};a_{\min})\big|^2 \\
	& \quad + 2 \big|V(Y_{\min};a)-V(Y_{\min};a_{\min})\big|^2 + 2 \big|V(Y_{\min};a_{\min})\big|^2 + \big(2\rho^{-\sigma}\big)^{(m+1)/(m-1)}.
\end{align*}
Combining~\eqref{contg} and the above inequality yields
\begin{align*}
	(V'(r;a))^2 & + \frac{\rho^\sigma}{2} V^2(r;a) \le \frac{2A_{\min}^2}{Y_{\min}^{N-1}} \big(a-a_{\min}\big)^2 + 2 \big|V'(Y_{\min};a_{\min})\big|^2 \\
	& \quad + 2 A_{\min}^2 \big(a-a_{\min}\big)^2 + 2 \big|V(Y_{\min};a_{\min})\big|^2 + \big(2\rho^{-\sigma}\big)^{(m+1)/(m-1)},
\end{align*}
hence
\begin{equation*}
	(V'(r;a))^2 + V^2(r;a) \le L(a_{\min},a_{\max},\rho), \quad r\in [Y_{\min},\rho]. 
\end{equation*}
We now combine the above inequality with~\eqref{eq1.rads} \blm{to complete the proof of Lemma~\ref{lem.gb}.}
\end{proof}

The next result establishes the precise continuity result with respect to the parameter $a$.
\begin{proposition}\label{prop.stab}
Let $a_0>0$. For any $\rho\in (0,S(a_0))$, there holds
\begin{equation}
	\lim_{a\to a_0} \sup_{r\in [0,\rho]} \left\{ |V(r;a)-V(r;a_0)| + r^{N-1} |V'(r;a)-V'(r,a_0)| \right\} = 0. \label{gs}
\end{equation}
In addition, if $S(a_0)<\infty$, then~\eqref{gs} is valid for $\rho=S(a_0)$.
\end{proposition}

\begin{proof}
We set $a_{\min}=a_0/2$, $a_{\max}=2a_0$ and consider a sequence $(a_j)_{j\ge 1}$ such that
\begin{equation}
	a_j\in [a_{\min},a_{\max}] \;\text{ for }\; j\ge 1 \;\;\text{ and }\;\; \lim_{j\to\infty} a_j=a_0. \label{b10}
\end{equation}
On the one hand, Lemma~\ref{lem.loccd} guarantees the existence of $Y_{\min}>0$ and $A_{\min}>0$ depending only on $a_{\min}$ such that
\begin{equation}\label{b15}
\begin{split}
	\sup\limits_{r\in[0,Y_{\min}]}|V(r;a_j)-V(r;a_0)| & \le A_{\min} |a_j-a_0|, \quad j\ge 1, \\
	\sup\limits_{r\in[0,Y_{\min}]} \left\{ r^{N-1} |V'(r;a_j)-V'(r;a_0)| \right\} & \le A_{\min} |a_j-a_0|, \quad j\ge 1.
\end{split}
\end{equation}
On the other hand, Lemma~\ref{lem.gb}, \blm{along with~\eqref{eq1.rads} and~\eqref{b15},} ensures that, for each $n\ge 2$, $(V(\cdot;a_j))_{j\ge 1}$ is a bounded sequence in $C^2([1/n,n])$. We then infer from the Ascoli-Arzel\`a Theorem and a diagonal process that there are $\bar{V}\in C^1((0,\infty))$ and a subsequence $(V(\cdot;a_{j_k}))_{k\ge 1}$ of $(V(\cdot;a_j))_{j\ge 1}$ such that
\begin{equation}
	\lim_{k\to\infty} \big|V(\cdot;a_{j_k}) - \bar{V}\big|_{C^1([1/n,n])} = 0, \quad n\ge 2. \label{b16}
\end{equation}
A first consequence of~\eqref{b16} is that 
\begin{equation}
	\lim_{k\to \infty} \sup_{r\in [1/n,n]} \left| V''(r;a_{j_k}) + \frac{N-1}{r} \bar{V}'(r) + r^\sigma \bar{V}(r) - \big(|\bar{V}|^{(1-m)/m}\bar{V}\big)(r) \right| = 0, \quad n\ge 2. \label{b17}
\end{equation}
We then deduce from~\eqref{b16} and~\eqref{b17} that $\bar{V}\in C^2((0,\infty))$ and is a solution \blm{to~\eqref{eq1s.rad}}. Furthermore, we infer from~\eqref{b15} that, for $r\in (0,Y_{\min})$ and $k\ge 1$,
\begin{align*}
	\big|V(r;a_{j_k})-a_0\big| & \le \big|V(r;a_{j_k})-V(r;a_0)\big| + \big|V(r;a_0)-a_0\big| \\
	& \le A_{\min} |a_{j_k}-a_0| + \big|V(r;a_0)-a_0\big|,
\end{align*}
and
\begin{align*}
	r^{N-1} \big|V'(r;a_{j_k})\big| & \le r^{N-1} \big|V'(r;a_{j_k}) - V(r;a_0)\big|  + r^{N-1} \big|V'(r;a_0\blm{)}\big| \\
	& \le A_{\min} |a_{j_k}-a_0| + r^{N-1} \big|V'(r;a_0)\big| .
\end{align*}
We then let $k\to\infty$ in the above two inequalities and use~\eqref{b10} and~\eqref{b16} to obtain
\begin{equation*}
	\big|\bar{V}(r)-a_0\big| \le \big|V(r;a_0)-a_0\big| \;\;\text{ and }\;\; r^{N-1} \big|\bar{V}'(r)\big| \le r^{N-1} \big|V'(r;a_0)\big|
\end{equation*}
for $r\in (0,Y_{\min})$. Recalling that $V(\cdot;a_0)$ satisfies~\eqref{icas}, we conclude that
\begin{equation}
	\lim_{r\to 0} \big|\bar{V}(r)-a_0\big| = \lim_{r\to 0} r^{N-1} \big|\bar{V}'(r)\big| = 0. \label{icb}
\end{equation}
Consequently, $\bar{V}$ is a solution \blm{to~\eqref{eq1s.rad}} with initial conditions~\eqref{icb} and it follows from Proposition~\ref{prop.CP} that $\bar{V}=V(\cdot;a_0)$ on $[0,S(a_0))$. Since $V(\cdot;a_0)$ is the only cluster point of $(V(\cdot;a_j))_{j\ge 1}$ in $C^1([1/n,\rho])$ for all $\rho\in (0,S(a_0))$ and $n>1$, classical compactness arguments entail that the whole sequence $(V(\cdot;a_j))_{j\ge 1}$ converges to $V(\cdot;a_0)$ in that space, and actually that 
\begin{equation*}
	\lim_{a\to a_0} \sup_{r\in [1/n,\rho]} \left\{ |V(r;a)-V(r;a_0)| + r^{N-1} |V'(r;a)-V'(r,a_0)| \right\} = 0
\end{equation*}
for all $\rho\in (0,S(a_0))$ and $n>1$. The extension of the above convergence to $[0,\rho]$ is provided by Lemma~\ref{lem.loccd}. Finally, if $S(a_0)<\infty$, then $\bar{V}=V(\cdot;a_0)$ on $[0,S(a_0)]$ according to Proposition~\ref{prop.CP} and we may choose $\rho=S(a_0)$ in the previous argument, thereby completing the proof.
\end{proof}
We are now in a position to introduce the variation function $\partial_a V(\cdot;a)$. Although it formally satisfies~\eqref{varfun}, at this point we do not know yet whether $a\mapsto V(\cdot;a)$ is differentiable with respect to $a$. In order to make rigorous the definition of $\partial_a V(\cdot;a)$, we introduce the solution $\varphi(\cdot;a)$ to the Cauchy problem~\eqref{CPphi} and then prove, a posteriori, that~\eqref{varfun} is satisfied from a rigorous point of view. 

We thus start with the well-posedness of the Cauchy problem~\eqref{CPphi} and establish the following result.
\begin{proposition}\label{prop.wpvphib}
Let $m$ and $\sigma$ be as in~\eqref{range.exp} and $a\in (0,\infty)$. Then, there is a unique solution $\varphi(\cdot;a)\in C([0,R(a))) \cap C^2((0,R(a)))$ to~\eqref{CPphi} and
\begin{equation}
	\varphi'(r;a) = - \frac{r^{\sigma+1}}{\sigma+N} + o(r^{\sigma+1}), \quad {\rm as} \ r\to 0. \label{phip0b}
\end{equation}
\end{proposition}

\begin{proof}
The Cauchy problem~\eqref{CPphi} is equivalent to the Cauchy problem~\eqref{singoper} with
\begin{equation*}
	U(r)=\frac{1}{m}V^{(1-m)/m}(r;a), \quad F\equiv0, \quad h_0=1, \quad r\in (0,R(a)).
\end{equation*}
Consequently, the existence and uniqueness of $\varphi(\cdot;a)$ are a direct  consequence of Lemma~\ref{lem.exist}. 
\end{proof}

The next result closes the circle by showing that the function $\varphi(\cdot;a)$, defined as above, is indeed the variation with respect to $a$ of $V(\cdot;a)$. Without the presence of a singular coefficient (such as $r^{\sigma}$) in Eq. \eqref{eq1.rad}, the following result would have directly followed by classical arguments. Since we have a coefficient that blows up as $r\to 0$, we have to provide supplementary work in order to prove the continuous (and then differentiable) dependence with respect to $a$ in a neighborhood of $r=0$, as it does not follow directly from standard theory.
\begin{proposition}\label{prop.phi}
The mapping
\begin{equation*}
	a\mapsto V(\cdot;a)
\end{equation*}
is differentiable with respect to $a\in(0,\infty)$ and we have
\begin{equation}\label{vardif}
 \partial_aV(r;a)=\varphi(r;a), \quad (r,a)\in(0,R(a))\times(0,\infty),
\end{equation}
recalling that $\varphi(\cdot;a)\in C([0,R(a)))\cap C^2((0,R(a)))$ is the unique solution to~\eqref{CPphi}.

Furthermore, $\varphi(\cdot;a)$ and $r\mapsto r^{N-1} \varphi'(r;a)$ depend continuously on the parameter $a>0$ on $[0,R(a))$.
\end{proposition}

\begin{proof}
We fix an arbitrary $a_{\min}\in (0,\infty)$ and, for $(a_0,a)\in (a_{\min},\infty)^2$, we set
\begin{equation*}
	\Xi(r):=\frac{V(r;a)-V(r;a_0)}{a-a_0}-\varphi(r;a_0), \quad r\in [0,Y_{\min}],
\end{equation*}
recalling that $Y_{\min}$ is defined in Lemma~\ref{lem.loccd}. It is then easy to derive from~\eqref{eq1.rad} and~\eqref{eq1.phi} that $\Xi$ is a solution to the Cauchy problem
\begin{subequations}\label{cont5}
\begin{equation}\label{cont5a}
\begin{split}
	\Xi''(r)+\frac{N-1}{r}\Xi'(r)&+r^{\sigma}\Xi(r)-D(r;a,a_0)\Xi(r)\\
	&= \left( D(r;a,a_0) - \frac{1}{m} V^{(1-m)/m}(r;a_0) \right) \varphi(r;a_0), \quad r\in [0,Y_{\min}],
\end{split}
\end{equation}
\begin{equation}\label{cont5b}
	\Xi(0)=0, \quad \lim\limits_{r\to0} r^{N-1}\Xi'(r)=0,
\end{equation}
\end{subequations}
with $D(\cdot;a,a_0)$ defined in~\eqref{D}. Owing to~\eqref{cont1}, \eqref{cont3} and the boundedness of $\varphi(\cdot;a_0)$ on $[0,Y_{\min}]$, we apply Lemma~\ref{lem.uniq} (with $U=D(\cdot;a,a_0)$ and $F=[D(\cdot;a,a_0)-V^{(1-m)/m}(\cdot;a_0)/m] \blm{\varphi(\cdot;a_0)}$) to~\eqref{cont5} to obtain
\begin{align*}
	|\Xi(r)|&\leq C_0(Y_{\min}) \left\| \left( D(\cdot;a,a_0) - \frac{1}{m} V^{(1-m)/m}(\cdot;a_0) \right) \varphi(\cdot;a_0) \right\|_{L^1(0,Y_{\min})}\\
	&\quad \times\exp\left\{ C_0(Y_{\min}) \left( \frac{Y_{\min}^{\sigma_N}}{\sigma_N} + Y_{\min}^{\sigma_N-\sigma-1} \|D(\cdot;a,a_0)\|_{L^1(0,Y_{\min})} \right) \right\}.
\end{align*}
Hence, in view of~\eqref{cont3},
\begin{equation}\label{cont6}
	|\Xi(r)| \leq B \left\| \left( D(\cdot;a,a_0) - \frac{1}{m} V^{(1-m)/m}(\cdot;a_0) \right) \varphi(\cdot;a_0) \right\|_{L^\infty(0,Y_{\min})},
\end{equation}
with
\begin{equation*}
	B := Y_{\min} C_0(Y_{\min}) \exp\left\{ C_0(Y_{\min}) \left( \frac{Y_{\min}^{\sigma_N}}{\sigma_N} + Y_{\min}^{\sigma_N-\sigma} \delta_{\min}^{(1-m)/m} \right) \right\}.
\end{equation*}
Using~\eqref{cont1}, we next estimate
\begin{equation*}
\begin{split}
	\left|V^{1/m}(r;a)-V^{1/m}(r;a_0)\right.&\left.-\frac{1}{m}(V(r;a)-V(r;a_0))V^{(1-m)/m}(r;a_0)\right|\\
	&= \frac{1}{m} \left|\int_{V(r;a_0)}^{V(r;a)}\left(z^{(1-m)/m}-V^{(1-m)/m}(r_0;a)\right)\,dz\right|\\
	&\leq\frac{m-1}{m^2}\left|\int_{V(r;a_0)}^{V(r;a)}\min\{z,V(r;a_0)\}^{(1-2m)/m}(z-V(r;a_0))\,dz\right|\\
	&\leq\frac{m-1}{2m^2}\delta_{\min}^{(1-2m)/m}|V(r;a)-V(r;a_0)|^2,
\end{split}
\end{equation*}
and the latter bound, along with~\eqref{cont4}, readily ensures, after dividing both sides by $\blm{|V(r;a)-V(r;a_0)|}$, that
\begin{equation*}
	\left\|D(\cdot;a,a_0)-\frac{1}{m}V^{(1-m)/m}(\cdot;a_0)\right\|_{L^{\infty}(0,Y_{\min})} \leq \delta_{\min}^{(1-2m)/m} A |a-a_0|.
\end{equation*}
Plugging the latter estimate in~\eqref{cont6} gives
\begin{equation}\label{cont7}
	\sup\limits_{r\in [0,Y_{\min}]} \left| \frac{V(r;a)-V(r;a_0)}{a-a_0} - \varphi(r;a_0) \right| \le B_{min} \|\varphi(\cdot;a_0)\|_{L^\infty(0,Y_{\min})} |a-a_0|
\end{equation}
for $(a_0,a)\in (a_{\min},\infty)^2$, with $B_{min} := AB \delta_{\min}^{(1-2m)/m}$, and~\eqref{vardif} follows immediately from~\eqref{cont7} on $[0,Y_{\min}]$.

We next set $\Psi := \varphi(\cdot;a) - \varphi(\cdot;a_0)$ and deduce from~\eqref{eq1.phi} that $\Psi$ solves the initial value problem
\begin{subequations}\label{CPPsi}
\begin{equation}\label{CPPsia}
\begin{split}
	\Psi''(r) + \frac{N-1}{r} \Psi'(r) + r^{\sigma} \Psi(r) & - \frac{V^{(1-m)/m}(r;a)}{m} \Psi(r) \\
	& = - \mathcal{D}(r;a,a_0) \varphi(r;a_0), \quad r\in (0,Y_{\min}),
\end{split}
\end{equation}
\begin{equation}\label{CPPsib}
	\Psi(0)=0, \quad \lim\limits_{r\to 0} r^{N-1} \Psi'(r)=0,
\end{equation}
\end{subequations}
with
\begin{equation*}
	\mathcal{D}(r;a,a_0) := \frac{1}{m} \left( V^{(1-m)/m}(r;a) - V^{(1-m)/m}(r;a_0) \right), \quad r\in [0,Y_{\min}].
\end{equation*}
We then argue as in the proof of Lemma~\ref{lem.loccd} and apply Lemma~\ref{lem.uniq} (with $R=Y_{\min}$, $U=V^{(1-m)/m}(\cdot;a)$ and $F=-\mathcal{D}(\cdot;a,a_0) \varphi(\cdot;a_0)$) to obtain that, in view of~\eqref{cont1}, 
\begin{equation*}
	|\Psi(r)| + r^{N-1} |\Psi'(r;a)| \le C_{\min} \|\varphi(\cdot;a_0)\|_{L^\infty(0,Y_{\min})} |a-a_0|, \quad r\in [0,Y_{\min}],
\end{equation*}
for some $C_{\min}>0$ depending only on $a_{\min}$. We have thus established the continuous dependence of $\varphi(\cdot;a)$ and $r\mapsto r^{N-1} \varphi'(r;a)$ with respect to $a$ on $[0,Y_{\min}]$. Since the function $r\mapsto r^\sigma$ has no singularity beyond $Y_{\min}$, we complete the proof by employing the standard results of continuity and differentiability with respect to the parameter in order to extend~\eqref{vardif} and the continuous dependence of $\varphi(\cdot;a)$ with respect to $a$ on $[Y_{\min},R(a_0))$.
\end{proof}

\bigskip

\section*{Acknowledgements} This work is partially supported by the Spanish project PID2024-160967NB-I00 funded by Agencia Estatal de Investigaci\'on (Spain) and FEDER. Part of this work has been developed during visits of R. G. I. to Laboratoire de Math\'ematiques LAMA, Universit\'e de Savoie \blm{Mont Blanc}, and of Ph. L. to Universidad Rey Juan Carlos, and both authors thank these institutions for hospitality and support.

For the purpose of Open Access, a CC-BY license has been applied by the authors to all versions of this article leading up to the Author Accepted Manuscript.

\bigskip

\noindent \textbf{Data availability} Our manuscript has no associated data.

\bigskip

\noindent \textbf{Conflict of interest} The authors declare that there is no conflict of interest.

\bibliographystyle{plain}



\end{document}